\documentclass[12pt,a4paper]{ctexart}  % 原 cctart(GBK) 已改为 ctexart(UTF-8), 适配现代 TeX Live

\usepackage{amssymb}
\usepackage[dvipsnames,usenames]{color}
\usepackage{amsfonts}
\usepackage{amssymb, amsmath}
\usepackage{latexsym,bm}
\usepackage{graphicx}
\usepackage{subfigure}
\usepackage{multirow}
\usepackage{color,xcolor} %\textcolor
\usepackage{amsthm} %包含证明环境
\usepackage[acronym]{glossaries}
\usepackage{float}

\newtheorem{theorem}{Theorem}[section]
\newtheorem{remark}[theorem]{Remark}
\newtheorem{lemma}[theorem]{Lemma}
\newtheorem{corollary}[theorem]{Corollary}

\newtheorem{assumption}{Assumption}[section]  % 不共同编号

\newacronym{DFT}{DFT}{discrete Fourier transform}
\newacronym{FFT}{FFT}{fast Fourier transform}
\newacronym{LOBPCG}{LOBPCG}{locally optimal block preconditioned conjugate gradient}
\newacronym{IPR}{IPR}{inverse participation ratio}

\begin{document}
	\begin{center}
		{\LARGE\textbf{Magnetic Bloch eigenvalue problems in nonuniform periodic magnetic fields: a Fourier–finite difference approach}}
	\end{center}
	\begin{center}
		{\large Xueyang Li\qquad Yuanwei Lai}\\[2pt]
		\texttt{lixy1217@xtu.edu.cn}\qquad\texttt{laiyuanwei02@163.com}
	\end{center}
	\begin{center}
		{\small Hunan Key Laboratory for Computation and Simulation in Science and Engineering, 
			Key Laboratory of Intelligent Computing and Information Processing of Ministry of Education, School of Mathematics and Computational Science,
			Xiangtan University, Xiangtan, Hunan 411105, P.~R.~China}
	\end{center}
	
	{\small\noindent\textbf{Abstract.} We study magnetic Bloch eigenvalue problems for two-dimensional magnetic Schr\"odinger operators with nonuniform periodic magnetic fields and periodic scalar potentials, and develop a hybrid Fourier–finite difference discretization. A Landau-type gauge is constructed by solving a periodic Poisson equation. Under the integer magnetic flux condition, the magnetic Bloch boundary conditions are incorporated through shifts of Fourier-mode indices. The method combines link-variable-based covariant central differences in the $x$-direction with a Fourier pseudospectral discretization in the $y$-direction, yielding a Hermitian discrete Hamiltonian. The resulting eigenvalue problems are solved using a \gls{LOBPCG} method. Under appropriate regularity assumptions, we establish a fully discrete consistency residual estimate of order $O(h_1^2+N_2^{-s})$. With suitable spectral separation and eigenspace dimension-preserving conditions, error estimates of the same order are derived for isolated eigenvalues and their associated eigenspaces. Numerical experiments show second-order convergence in the finite difference direction and rapid error reduction with increasing Fourier resolution, and illustrate band degeneracies and eigenstate distributions in nonuniform periodic magnetic fields. \par}
	
	\glsresetall
	
	\vspace{2mm}
	{\small\noindent\textbf{Mathematics Subject Classification (2020):} 65N25, 65N35, 65N15, 35P15.\par}
	
	\vspace{2mm}
	{\small\noindent\textbf{Keywords:} magnetic Schr\"odinger eigenvalue problem; magnetic Bloch boundary conditions; link-variable-based covariant central differences; Fourier pseudospectral method; eigenvalue and eigenspace error estimates.\par}
	\vskip 0.8cm
	
	\section{Introduction}\setcounter{equation}{0}
	
	Magnetic Schr\"odinger operators describe the motion of charged quantum particles and play an important role in condensed matter physics, superconductivity, the quantum Hall effect, and electronic structure calculations. In a periodic
	lattice potential, an applied magnetic field alters both the local oscillations and the spatial
	structure of eigenstates; the compatibility between magnetic flux and lattice periodicity also
	shapes the band structure. Hofstadter's study of two-dimensional periodic systems at rational and
	irrational magnetic flux revealed the well-known fractal structure of the
	spectrum \cite{Hofstadter1976}. Thouless, Kohmoto, Nightingale, and den Nijs related the Hall
	conductance of two-dimensional periodic systems to topological invariants of the energy
	bands \cite{ThoulessKohmotoNightingaleDenNijs1982}. Chern numbers can be computed efficiently on a
	discretized Brillouin zone using a gauge-invariant scheme \cite{FukuiHatsugaiSuzuki2005}. The
	spatial localization of electronic states in magnetic fields remains an important problem in the quantum
	Hall effect \cite{Seye2025}. In the study of band dynamics, Nenciu gave a systematic, rigorous
	treatment of effective Hamiltonians for Bloch electrons in electromagnetic fields \cite{Nenciu1991},
	while Panati, Spohn, and Teufel investigated the effective dynamics of Bloch electrons and the
	mathematical foundations of the Peierls substitution \cite{PanatiSpohnTeufel2003}.
	
	For a magnetic field with nonzero mean, lattice translations must be accompanied by phase factors
	determined by the vector potential. This distinguishes the problem from the usual periodic
	Schr\"odinger setting and requires magnetic translations and magnetic Bloch theory in place of
	conventional Bloch theory. Brown's work on Bloch electrons in uniform magnetic
	fields \cite{Brown1964} and Zak's theory of the magnetic translation group \cite{Zak1964a} laid the
	foundations of this framework. Gruber examined Bloch theory for magnetic systems from a broader
	geometric and spectral perspective \cite{Gruber2000}. In spectral theory, Avron and Simon studied
	the stability of spectral gaps under variations in the magnetic field in the presence of a periodic
	potential \cite{AvronSimon1985}. Helffer and Sj\"ostrand established connections between magnetic
	Schr\"odinger operators and the Harper equation \cite{HelfferSjostrand1989}, and Cornean
	investigated the continuity of spectral bands of Harper-type and magnetic Schr\"odinger operators
	with respect to the magnetic field parameter \cite{Cornean2010}. 
	
	Several classes of numerical methods have been developed for Schr\"odinger eigenvalue problems in
	electromagnetic fields. Alouges and Bonnaillie-No\"el used finite element methods to study
	ground-state computation and gauge-invariant discretization for constant magnetic
	fields \cite{AlougesBonnaillieNoel2006}. Bonnaillie-No\"el, Dauge, Martin, and Vial performed
	high-order finite element computations of low-energy eigenstates on domains with corners, focusing
	on eigenfunction localization and high-frequency oscillations in strong magnetic
	fields \cite{BonnaillieNoelDaugeMartinVial2007}. These studies concern bounded-domain problems with
	constant magnetic fields. Another class of methods introduces
	link variables to preserve gauge invariance at the discrete level. Governale and Ungarelli
	constructed a gauge-invariant grid discretization \cite{GovernaleUngarelli1998}. Christiansen and
	Halvorsen combined ideas from lattice gauge theory with finite elements on simplicial meshes to
	construct a gauge-invariant discretization of the Schr\"odinger eigenvalue problem in an
	electromagnetic field and analyze its error \cite{Christiansen}. Halvorsen and Kvaal subsequently
	considered manifestly gauge-invariant constructions applicable to general grids and even Fourier
	pseudospectral discretizations \cite{HalvorsenKvaal2012}. These studies provide a basis for
	structure-preserving discretizations using link variables. Cai and Galli introduced a periodic-supercell
	method for first-principles calculations in uniform magnetic fields \cite{CaiGalli2004}, and Lee,
	Cai, and Galli subsequently developed the corresponding electronic structure
	algorithms \cite{LeeCaiGalli2007}. These methods exploit Fourier representations for uniform magnetic fields, but extending them to spatially varying periodic fields while retaining their spectral accuracy and computational efficiency is not straightforward. More recently, Ovall and Zhu used
	gauge invariance to construct a canonical gauge for eigenvalue computations by solving a Poisson
	problem, with the aim of reducing eigenfunction oscillations and improving computational efficiency
	and stability \cite{OvallZhu2025}.
	
	We develop a hybrid Fourier--finite difference discretization for two-dimensional magnetic Schr\"odinger eigenvalue problems with nonuniform periodic magnetic fields, periodic scalar potentials, and magnetic Bloch boundary conditions. The scheme treats the variable-coefficient covariant derivatives together with the magnetic phase factors at the boundaries. The method combines link-variable-based covariant central differences
	in the $x$-direction with a Fourier pseudospectral
	discretization in the $y$-direction. The
	resulting hybrid scheme preserves the Hermitian structure of the discrete Hamiltonian matrix and
	treats variable-coefficient terms efficiently using local finite-difference operations and
	one-dimensional \gls{FFT} in $y$. We also derive error estimates for isolated
	eigenvalues and their eigenspaces under suitable regularity assumptions, provided that the
	corresponding discrete spectral subspace has the same dimension as the target eigenspace.
	
	The remainder of the paper is organized as follows. Section 2 introduces a Landau-type gauge for
	nonuniform periodic magnetic fields, states the magnetic Bloch boundary conditions, and
	formulates the corresponding eigenvalue problem on a bounded magnetic unit cell. Section 3
	constructs the fully discrete scheme by combining covariant central differences in $x$ with a Fourier pseudospectral discretization in $y$, and discusses the
	construction of a preconditioner for the iterative eigensolver. Section 4 establishes estimates for the consistency residual and for
	the eigenvalue and eigenspace errors of the fully discrete scheme. Section 5 presents numerical experiments illustrating the spectral properties and convergence behavior of the proposed method.
	\vskip 0.5cm

	\section{Magnetic Bloch formulation and eigenvalue problem}
	\setcounter{equation}{0}
	
	Let $B$ and $V$ be real-valued doubly periodic functions representing the magnetic field and scalar potential, respectively. Specifically, for any $x,y\in\mathbb{R}$, we assume that $B(x,y),V(x,y)\in\mathbb{R}$ and that there exist positive constants $a,b$ such that
	\begin{equation*}
		\begin{aligned}
			&B(x+a,y)=B(x,y+b)=B(x,y),\\
			&V(x+a,y)=V(x,y+b)=V(x,y).
		\end{aligned}
	\end{equation*}
	Define the magnetic unit cell by
	\begin{equation*}
		\Omega:=\left(-\frac a2,\frac a2\right)\times\left(-\frac b2,\frac b2\right),
	\end{equation*}
	and denote the mean magnetic field by
	\begin{equation}\label{eq:mean B}
		\bar B:=\frac1{ab}\int_\Omega B(x,y)\,\mathrm{d}x\mathrm{d}y.
	\end{equation}
	Throughout this paper, we assume that the normalized magnetic flux through a magnetic unit cell satisfies the integer magnetic flux condition
	\begin{equation}\label{eq:integer flux}
		n_B:=\frac{\bar B ab}{2\pi}\in\mathbb Z.
	\end{equation}
	Here, $\bar B ab$ is the total magnetic flux through the unit cell, while $n_B$ denotes the normalized magnetic flux.
	If the magnetic flux through a primitive cell, measured in units of $2\pi$, is rational, the magnetic unit cell can be enlarged to satisfy \eqref{eq:integer flux}. For irrational magnetic flux, rational approximations obtained through Diophantine approximation can be used as substitutes.
	
	We next construct a vector potential associated with the magnetic field. Let $\phi$ be determined by the periodic Poisson problem
	\begin{equation}\label{eq:poisson phi}
		\left\{\begin{aligned}
			&\Delta\phi=B(x,y)-\bar B,\\
			&\phi(x+a,y)=\phi(x,y+b)=\phi(x,y),\\
			&\int_\Omega\phi(x,y)\,\mathrm{d}x\mathrm{d}y=0.
		\end{aligned}\right.
	\end{equation}
	We adopt the following Landau-type gauge:
	\begin{equation}\label{eq:vector potential}
		\bm A(x,y):=\big(-\phi_y(x,y),\ \bar Bx+\phi_x(x,y)\big),
	\end{equation}
	where $\phi_x=\partial_x\phi$ and $\phi_y=\partial_y\phi$. It follows that
	\begin{equation}\label{eq:curl}
		\partial_xA_y-\partial_yA_x=\bar B+\Delta\phi=B(x,y).
	\end{equation}
	In the discretization and error analysis below, $\bar B$ is treated as a given real parameter, and $\phi_x$ and $\phi_y$ as known functions satisfying the stated regularity assumptions. Errors arising from the numerical solution of \eqref{eq:poisson phi} are not included.
	
	The corresponding magnetic Schr\"odinger operator in the Landau-type gauge is given by
	\begin{equation}\label{eq:original H}
		H\psi:=\frac1{2}\left[
		\left(-\mathrm i\partial_x+\phi_y(x,y)\right)^2+
		\left(-\mathrm i\partial_y-\bar Bx-\phi_x(x,y)\right)^2
		\right]\psi+V\psi.
	\end{equation}
	For simplicity, physical constants such as Planck's constant, the particle charge, and the particle mass have been absorbed into the corresponding coefficients. The operator \eqref{eq:original H} is defined on an unbounded domain, and its spectral problem cannot be directly identified with a bounded-domain eigenvalue problem having a discrete spectrum. According to the theory of magnetic translation groups and magnetic Bloch decomposition \cite{Brown1964,Zak1964a,Gruber2000}, the spectral problem on the unbounded domain can be reduced to a family of parameter-dependent eigenvalue problems on the magnetic unit cell. For a given quasimomentum $\bm{k}=(k_1,k_2)$, the original wave function can be expressed as
	\begin{equation}\label{eq:Bloch peel}
		\psi(x,y)=\mathrm e^{\mathrm i(k_1x+k_2y)}u(x,y).
	\end{equation}
	This transformation leads to the parameter-dependent operator acting on $u$:
	\begin{equation}\label{eq:Hk continuous}
		H^{(\bm{k})}u:=\frac1{2}\left(D_x^2+D_y^2\right)u+V(x,y)u,
	\end{equation}
	where
	\begin{equation}\label{eq:Dx continuous}
		D_x:=-\mathrm i\partial_x+k_1+\phi_y(x,y)
	\end{equation}
	and
	\begin{equation}\label{eq:Dy continuous}
		D_y:=-\mathrm i\partial_y+k_2-\bar Bx-\phi_x(x,y)
	\end{equation}
	are the first-order covariant differential operators in the $x$- and $y$-directions, respectively. Meanwhile, the boundary conditions for $u(x,y)$ take the following magnetic Bloch form:
	\begin{equation}\label{eq:magnetic bc}
		u(x+a,y)=\mathrm e^{\mathrm i\bar Bay}u(x,y),\qquad
		u(x,y+b)=u(x,y).
	\end{equation}
	The integer magnetic flux condition \eqref{eq:integer flux} implies that
	$\mathrm e^{\mathrm i\bar Ba(y+b)}=\mathrm e^{\mathrm i\bar Bay}$,
	so that the boundary conditions \eqref{eq:magnetic bc} are compatible with periodicity in the $y$-direction. Under these boundary conditions, the operator \eqref{eq:Hk continuous} defines a family of eigenvalue problems parameterized by $\bm{k}$. The corresponding magnetic Bloch eigenstates are recovered by multiplying the eigenfunctions by the phase factor in \eqref{eq:Bloch peel}. Collectively, these states provide a complete spectral representation of the original operator \eqref{eq:original H} in the direct-integral sense.
	
	We therefore consider the following eigenvalue problem on the bounded domain $\Omega$:
	\begin{equation}\label{eq:nonconst bounded eigenproblem}
		\left\{\begin{aligned}
			&H^{(\bm{k})}u=\lambda^{(\bm{k})} u,\qquad (x,y)\in \Omega,\\
			&u\Big(\frac a2,y\Big)=\mathrm e^{\mathrm i\bar Bay}u\Big(-\frac a2,y\Big),\\
			&D_xu\Big(\frac a2,y\Big)=\mathrm e^{\mathrm i\bar Bay}D_xu\Big(-\frac a2,y\Big),\\
			&u\Big(x,\frac b2\Big)=u\Big(x,-\frac b2\Big),\\
			&D_yu\Big(x,\frac b2\Big)=D_yu\Big(x,-\frac b2\Big),\\
			&\|u\|_{L^2(\Omega)}=1,
		\end{aligned}\right.
	\end{equation}
	where $\|u\|_{L^2(\Omega)}:=\left<u,u\right>^{\frac12}$, with the inner product defined by
	\begin{equation*}
		\left<u,w\right>:=\int_\Omega u(x,y)\overline{w(x,y)}\,\mathrm{d}x\mathrm{d}y.
	\end{equation*}
	Under suitable regularity assumptions on the coefficients, the above magnetic Bloch boundary value problem defines a self-adjoint operator $H^{(\bm{k})}$ with compact resolvent. Consequently, its spectrum consists of eigenvalues of finite multiplicity, with $+\infty$ as the only possible accumulation point.
	
	Using the periodicity of $u(x,y)$ in the $y$-direction, we obtain the Fourier expansion
	\begin{equation}\label{eq:y Fourier continuous}
		u(x,y)=\sum_{j\in\mathbb Z}\widehat u_j(x)\mathrm e^{\mathrm i\kappa_jy},
		\qquad
		\kappa_j=\frac{2\pi j}{b},
	\end{equation}
	where the Fourier coefficients in the $y$-direction are given by
	\begin{equation*}
		\widehat u_j(x):=\frac1b\int_{-b/2}^{b/2}u(x,y)\mathrm e^{-\mathrm i\kappa_jy}\,\mathrm{d}y.
	\end{equation*}
	By the integer magnetic flux condition \eqref{eq:integer flux}, we have
	\[
	\mathrm e^{\mathrm i\bar Bay}=\mathrm e^{2\pi\mathrm i n_By/b}.
	\]
	Combining this identity with \eqref{eq:magnetic bc} yields the magnetic Bloch boundary condition in terms of Fourier coefficients:
	\begin{equation}\label{eq:coeff magnetic bc}
		\widehat u_j\left(\frac a2\right)
		=\widehat u_{j-n_B}\left(-\frac a2\right),
		\qquad j\in\mathbb Z.
	\end{equation}
	Similarly, the magnetic Bloch boundary condition for the first-order covariant derivative in the $x$-direction in \eqref{eq:nonconst bounded eigenproblem} can be expressed as
	$$
	\widehat{D_xu}_j\left(\frac a2\right)
	=\widehat{D_xu}_{j-n_B}\left(-\frac a2\right),
	\qquad j\in\mathbb Z.
	$$

	\section{Fourier--Finite Difference Discretization}\setcounter{equation}{0}
	
	\subsection{Grid functions and discrete Fourier coefficient vectors}
	Let $N_1$ and $N_2$ be positive even integers. We introduce the grid points
	\begin{equation}\label{eq:grid nodes}
		\begin{aligned}
			&x_l:=-\frac a2+lh_1,\qquad l=0,1,\ldots,N_1-1,\\
			&y_m:=-\frac b2+mh_2,\qquad m=0,1,\ldots,N_2-1,
		\end{aligned}
	\end{equation}
	with mesh sizes
	\begin{equation*}
		\bm h=(h_1,h_2):=\Big(\frac a{N_1},\frac b{N_2}\Big).
	\end{equation*}
	A vector of values at these grid points,
	$$\bm U:=(U_{l,m})_{0\leq l\leq N_1-1,~0\leq m\leq N_2-1}$$
	is referred to as a grid function. For each fixed $l$, we also write
	$$\bm U_l=(U_{l,m})_{0\leq m\leq N_2-1}.$$
	We use the same grid-point definitions \eqref{eq:grid nodes} for $l,m\in\mathbb Z$ when indices outside the computational grid are required.
	
	Let
	\begin{equation*}
		\mathcal K_{N_2}:=\left\{-\frac{N_2}{2},\,\ldots\,,\,\frac{N_2}{2}-1\right\}
	\end{equation*}
	denote the set of Fourier-mode indices in the $y$-direction.
	For each fixed $x_l$, the grid values $\bm U_l$ are represented by the
	following discrete Fourier coefficient vector in the $y$-direction:
	$$\widetilde {\bm U}_l=(\widetilde U_{l,j})_{j\in\mathcal K_{N_2}}=\mathbf{F}_y\bm U_l,$$
	where the \gls{DFT} in the $y$-direction is defined by
	\begin{equation}\label{eq:DFT}
		\big(\mathbf{F}_y\bm U_l\big)_j:=\frac1{N_2}\sum_{m=0}^{N_2-1}
		U_{l,m}\mathrm e^{-\mathrm i\kappa_jy_m},
		\qquad j\in\mathcal K_{N_2},
	\end{equation}
	and its inverse is given by
	\begin{equation*}
		\big(\mathbf{F}_y^{-1}\widetilde {\bm U}_l\big)_m=\sum_{j\in\mathcal K_{N_2}}
		\widetilde U_{l,j}\mathrm e^{\mathrm i\kappa_jy_m},
		\qquad m=0,1,\ldots,N_2-1.
	\end{equation*}
	The definition of the \gls{DFT} immediately yields the discrete Parseval identity
	\begin{equation}\label{eq:parseval}
		h_2\sum_{m=0}^{N_2-1}U_{l,m}\overline{W_{l,m}}
		=b\sum_{j\in\mathcal K_{N_2}}
		\widetilde U_{l,j}\overline{\widetilde W_{l,j}}.
	\end{equation}
	For subsequent use, we collect the discrete Fourier coefficients into the vector
	\begin{equation}\label{eq:U vector}
		\widetilde{\bm U}=\big(\widetilde U_{l,j}\big)_{0\le l\le N_1-1,\ j\in\mathcal K_{N_2}}.
	\end{equation}
	We equip the space of discrete Fourier coefficient vectors with the inner product
	\begin{equation*}
		\big<\widetilde{\bm U},\widetilde{\bm W}\big>_{\bm h}
		:=h_1b\sum_{l=0}^{N_1-1}\sum_{j\in\mathcal K_{N_2}}
		\widetilde U_{l,j}\overline{\widetilde W_{l,j}},
	\end{equation*}
	and the associated norm
	\begin{equation*}
		\big\|\widetilde{\bm U}\big\|_{\bm h}
		:=\big<\widetilde{\bm U},\widetilde{\bm U}\big>_{\bm h}^{1/2}.
	\end{equation*}
	For a function $\bm g=(g_m)_{m=0}^{N_2-1}$ specified on the $y$-grid $\{y_m\}_{m=0}^{N_2-1}$, its multiplication operator in the discrete Fourier coefficient space is defined by
	\begin{equation}\label{eq:multiplication operator}
		\mathbf M_{\bm g}:=\mathbf F_y\operatorname{diag}(g_0,g_1,\ldots,g_{N_2-1})\mathbf F_y^{-1}.
	\end{equation}
	If all components of $\bm g$ are real, then the matrix $\mathbf M_{\bm g}$ defined by \eqref{eq:multiplication operator} is Hermitian with respect to the Fourier inner product on each fixed $x$-layer.
	
	We now transfer the magnetic Bloch boundary conditions \eqref{eq:coeff magnetic bc}, expressed in continuous Fourier coefficients, to the discrete Fourier coefficient space. By the periodicity of the \gls{DFT}, the definition \eqref{eq:DFT} extends naturally to indices $j\notin\mathcal K_{N_2}$ and satisfies
	$\widetilde U_{l,j}=\widetilde U_{l,\,(j~\mathrm{MOD}~N_2)}$. Consequently,
	the discrete magnetic Bloch boundary conditions take the form
	\begin{equation}\label{eq:ghost Fourier}
		\widetilde{U}_{N_1,j}=\widetilde{U}_{0,\,[j-n_B]_{\mathcal K_{N_2}}},\quad
		\widetilde{U}_{-1,j}=\widetilde{U}_{N_1-1,\,[j+n_B]_{\mathcal K_{N_2}}},\qquad\forall j\in\mathcal K_{N_2},
	\end{equation}
	Here, $[m]_{\mathcal K_{N_2}}$ denotes the unique integer in $\mathcal K_{N_2}$ congruent to $m$ modulo $N_2$. The quantities $\widetilde{U}_{N_1,j}$ and $\widetilde{U}_{-1,j}$ are not actual components of the vector \eqref{eq:U vector}; they are introduced solely as ghost values for expressing the boundary difference stencils. Equation \eqref{eq:ghost Fourier} gives the corresponding magnetic Bloch boundary conditions on the physical grid:
	\begin{equation*}
		U_{N_1,m}=\mathrm e^{\mathrm i\bar Bay_m}U_{0,m},
		\qquad
		U_{-1,m}=\mathrm e^{-\mathrm i\bar Bay_m}U_{N_1-1,m}.
	\end{equation*}
	Thus, the discrete extension agrees with the continuous magnetic Bloch boundary conditions \eqref{eq:magnetic bc} at the grid points.

	\subsection{Covariant central differences in $x$ and Fourier pseudospectral discretization in $y$}
	We discretize the second-order covariant derivative in the $x$-direction using link-variable-based covariant central differences, hereafter referred to as the covariant central difference scheme. The essential idea is to construct the covariant differences through phase parallel transport between neighboring grid points. Related approaches can be found in gauge-invariant lattice discretizations and link-variable-based structure-preserving methods \cite{Christiansen,GovernaleUngarelli1998,HalvorsenKvaal2012}. Define
	\begin{equation}\label{eq:theta integral exact}
		\vartheta_{l,m}(x):=\int_{x_l}^{x}\big(k_1+\phi_y(s,y_m)\big)\,\mathrm{d}s.
	\end{equation}
	A direct calculation shows that the covariant differential operator $D_x$ defined in \eqref{eq:Dx continuous}, together with its higher powers, can be represented by conjugating ordinary differential operators with phase factors:
	\begin{equation*}
		D_x^n u(x,y_m)=\mathrm e^{-\mathrm i\vartheta_{l,m}(x)}(-\mathrm i\partial_x)^n\big(\mathrm e^{\mathrm i\vartheta_{l,m}(x)}u(x,y_m)\big),\qquad n=1,2,\cdots.
	\end{equation*}
	Applying the trapezoidal rule to the integral in \eqref{eq:theta integral exact} gives
	\begin{equation*}
		\begin{aligned}
			\vartheta_{l,m}(x_{l+1})\approx\Theta_{l+1/2,m}
			&:=
			h_1k_1+\frac{h_1}{2}\big(\phi_y(x_l,y_m)+\phi_y(x_{l+1},y_m)\big),\\
			\mathrm e^{\mathrm i\vartheta_{l,m}(x_{l+1})}\approx G_{l+1/2,m}
			&:=
			\mathrm e^{\mathrm i\Theta_{l+1/2,m}}.
		\end{aligned}
	\end{equation*}
	Here, $G_{l+1/2,m}$ is the link variable associated with two neighboring grid points. Similarly,
	\begin{equation*}
		\mathrm e^{\mathrm i\vartheta_{l,m}(x_{l-1})}=\mathrm e^{-\mathrm i\vartheta_{l-1,m}(x_{l})}\approx \overline{G_{l-1/2,m}},
	\end{equation*}
	Set $\bm G_{l+1/2}=(G_{l+1/2,m})_{m=0}^{N_2-1}$ and $\bm G_{l-1/2}=(G_{l-1/2,m})_{m=0}^{N_2-1}$.
	These relations lead to the following covariant forward difference approximation of $D_x$, acting on Fourier coefficient vectors:
	$$\big(\mathbf{D}_{x,h_1}\widetilde{\bm U}\big)_{l}:=-\frac{\mathrm i}{h_1}\big(\mathbf M_{\bm G_{l+1/2}}\widetilde{\bm U}_{l+1}-\widetilde{\bm U}_l\big),\qquad l=0,1,\cdots, N_1-1.$$
	Its Hermitian adjoint is given by
	$$\big(\mathbf{D}_{x,h_1}^*\widetilde{\bm U}\big)_{l}:=\frac{\mathrm i}{h_1}\big(\mathbf M_{\bm G_{l-1/2}}^*\widetilde{\bm U}_{l-1}-\widetilde{\bm U}_l\big),\qquad l=0,1,\cdots, N_1-1.$$
	Under the discrete magnetic Bloch boundary conditions \eqref{eq:ghost Fourier}, the ghost vectors $\widetilde{\bm U}_{N_1}$ and $\widetilde{\bm U}_{-1}$ are uniquely determined by the components of \eqref{eq:U vector}. The difference operators are therefore well defined at $l=0$ and $l=N_1-1$ as well. The resulting covariant central difference approximation of $D_x^2$ is
	\begin{equation}\label{eq:Lx explicit interior}
		\begin{aligned}
			&(\mathbf L_{x,h_1}\widetilde{\bm U})_l:=\big(\mathbf{D}_{x,h_1}\mathbf{D}_{x,h_1}^*\widetilde{\bm U}\big)_{l}
			=\frac{2\widetilde{\bm U}_l-
				\mathbf M_{\bm G_{l+1/2}}\widetilde{\bm U}_{l+1}
				-\mathbf M_{\bm G_{l-1/2}}^*\widetilde{\bm U}_{l-1}}{h_1^2},\\
			&l=0,1,\cdots, N_1-1.
		\end{aligned}
	\end{equation}

	For each $x_l$, define
	\begin{equation*}
		\big(\bm A_l^y\big)_m:=k_2-\bar Bx_l-(\phi_x)_{l,m}.
	\end{equation*}
	The Fourier differentiation matrix in the $y$-direction is
	\begin{equation*}
		\bm{\Lambda}_y:=\operatorname{diag}(\kappa_j)_{j\in\mathcal K_{N_2}}.
	\end{equation*}
	The Fourier pseudospectral discretization of the second-order covariant differential operator $D_y^2$ is then given by
	\begin{equation}\label{eq:Ly square}
		\big(\mathbf L_{y,h_2}\widetilde{\bm U}\big)_l
		:=\big(\bm\Lambda_y+\mathbf M_{\bm A_l^y}\big)^2\widetilde {\bm U}_l,\qquad l=0,1,\cdots,N_1-1.
	\end{equation}
	The associated quadratic form satisfies
	\begin{equation*}
		\big<\mathbf L_{y,h_2}\widetilde{\bm U},\widetilde{\bm U}\big>_{\bm h}
		=h_1b\sum_{l=0}^{N_1-1}
		\big\|\big(\bm\Lambda_y+\mathbf M_{\bm A_l^y}\big)\widetilde{\bm U}_l\big\|^2\ge0.
	\end{equation*}

	The discrete potential operator, acting on Fourier coefficient vectors, is defined by
	\begin{equation*}
		\begin{aligned}
			\big(\mathbf{\widetilde{V}}_{\bm h}\widetilde{\bm U}\big)_l
			:=\mathbf{M}_{\bm V_l}\widetilde{\bm U}_l,\qquad l=0,1,\cdots, N_1-1,
		\end{aligned}
	\end{equation*}
	where $\bm V_l=(V(x_l,y_m))_{m=0}^{N_2-1}$. Combining \eqref{eq:Lx explicit interior} and \eqref{eq:Ly square}, we obtain the discrete Hamiltonian corresponding to the parameter-dependent magnetic Schr\"odinger operator \eqref{eq:Hk continuous}:
	\begin{equation}\label{eq:Hamilton discrete}
		\begin{aligned}
			&\mathbf H^{(\bm{k})}_{\bm h}\widetilde{\bm U}
			:=\frac1{2}\left(\mathbf L_{x,h_1}\widetilde{\bm U}+\mathbf L_{y,h_2}\widetilde{\bm U}\right)+\mathbf{\widetilde{V}}_{\bm h}\widetilde{\bm U}.
		\end{aligned}
	\end{equation}
	It follows from \eqref{eq:Lx explicit interior} and \eqref{eq:Ly square} that $\mathbf H^{(\bm{k})}_{\bm h}$ is Hermitian. Moreover, if $V(x,y)\geq0$, then $\mathbf H^{(\bm{k})}_{\bm h}$ is positive semidefinite.
	The discrete counterpart of the eigenvalue problem \eqref{eq:nonconst bounded eigenproblem} is therefore to find $\lambda_{\bm h}^{(\bm{k})}\in\mathbb R$ and a nonzero eigenvector $\widetilde{\bm U}_{\bm h}\in\mathbb C^{N_1N_2}$ such that
	\begin{equation}\label{eq:eigenproblem discrete}
		\mathbf H_{\bm h}^{(\bm{k})}\widetilde{\bm U}_{\bm h}
		=\lambda_{\bm h}^{(\bm{k})}\widetilde{\bm U}_{\bm h}.
	\end{equation}

	\subsection{Construction of the preconditioner}
	We solve the discrete eigenvalue problem \eqref{eq:eigenproblem discrete} iteratively using the \gls{LOBPCG} method. Introduced by Knyazev \cite{Knyazev2001LOBPCG}, this block eigensolver approximates several eigenvectors simultaneously. At each iteration, the local trial subspace is augmented with preconditioned residuals and search directions from the preceding iteration, and the approximations are updated through a Rayleigh--Ritz projection. As the iteration proceeds, the basis vectors of the trial subspace may become nearly linearly dependent, making the reduced Rayleigh--Ritz eigenvalue problem ill-conditioned. Orthogonalization of the basis vectors improves numerical stability. We employ the orthogonalized version of \gls{LOBPCG} described in Algorithm~2 of \cite{Hetmaniuk}.
	
	We seek a positive definite preconditioner $\mathbf P$ whose inverse can be applied efficiently to improve the residual-based search directions and accelerate convergence toward the desired eigenspace. When the variations in the scalar potential and periodic magnetic field are moderate and have a limited effect on conditioning, the variable coefficients can be approximated by constants in the preconditioning model. Accordingly, neglecting $\phi_x$ and $\phi_y$ in \eqref{eq:Dx continuous} and \eqref{eq:Dy continuous} leads to the following simplified counterparts of \eqref{eq:Lx explicit interior} and \eqref{eq:Ly square}, respectively:
	\begin{equation*}
		(\mathbf L_{x,h_1}^0\widetilde{\bm U})_l
		:=\frac{2\widetilde{\bm U}_l-\mathrm e^{\mathrm ik_1h_1}\widetilde{\bm U}_{l+1}
			-\mathrm e^{-\mathrm ik_1h_1}\widetilde{\bm U}_{l-1}}{h_1^2},
		\qquad 0\le l\le N_1-1,
	\end{equation*}
	and
	\begin{equation*}
		(\mathbf L_{y,h_2}^0\widetilde{\bm U})_{l}
		:=\left(\bm\Lambda_y+k_2\mathbf{I}_{N_2}-\bar Bx_l\mathbf{I}_{N_2}\right)^2\widetilde{\bm U}_{l},
	\end{equation*}
	where $\mathbf{I}_{N_2}$ denotes the identity operator on $\mathbb C^{N_2}$. We then define the preconditioner by
	\begin{equation}\label{eq:preconditioner P}
		\mathbf P:=\frac1{2}\left(\mathbf L_{x,h_1}^0+\mathbf L_{y,h_2}^0\right)
		+(\bar V+\sigma) \mathbf{I}_{N_1N_2},
	\end{equation}
	where $\bar V$ is the discrete average of the scalar potential,
	\begin{equation*}
		\bar V:=\frac1{N_1N_2}\sum_{l=0}^{N_1-1}\sum_{m=0}^{N_2-1}V_{l,m},
	\end{equation*}
	and the shift parameter $\sigma$ is chosen to ensure that $\mathbf P$ is positive definite. In particular, we generally impose $\sigma>-\bar V$.

	To characterize the sparsity pattern of $\mathbf P$ and develop an efficient inverse application, we first specify an ordering of the two-dimensional indices.
	Let $d=\gcd(|n_B|,N_2)$. The set $\{0,1,\cdots,N_2-1\}$ can then be partitioned into $d$ ordered sequences, each of length $N_2/d$:
	$$\mathcal{J}_n:=\Big(\iota\Big(\frac{N_2}{d}n\Big),\iota\Big(\frac{N_2}{d}n+1\Big),\cdots,\iota\Big(\frac{N_2}{d}n+\frac{N_2}{d}-1\Big)\Big),\quad n=0,1,\cdots,d-1,$$
	where $\iota\big(\frac{N_2}{d}n\big)=n$ and
	\begin{equation}\label{eq:index relation}
		\begin{aligned}
			&\iota(j)=\left[\iota(j-1)-n_B\right]_{\mathcal K_{N_2}},\qquad &\mathrm{for}~\frac{jd}{N_2}\notin\mathbb{Z},
		\end{aligned}
	\end{equation}
	Here, $[m]_{\mathcal K_{N_2}}$ denotes the representative of $m$ modulo $N_2$ in $\mathcal K_{N_2}$; that is, $[m]_{\mathcal K_{N_2}}\in\mathcal K_{N_2}$ and
	$$[m]_{\mathcal K_{N_2}}\,\mathrm{MOD}\,N_2=m\,\mathrm{MOD}\,N_2.$$
	The preceding relation further implies
	$$\begin{aligned}
		&\iota\Big(\frac{N_2}{d}n\Big)=\left[\iota\Big(\frac{N_2}{d}n+\frac{N_2}{d}-1\Big)-n_B\right]_{\mathcal K_{N_2}},
	\end{aligned}$$
	Thus, each $\mathcal{J}_n$
	forms a cyclic sequence under the relation \eqref{eq:index relation}.
	Concatenating these sequences gives the permutation
	$$\mathcal{J}:=(\mathcal{J}_n)_{n=0}^{d-1}=(\iota(m))_{m=0}^{N_2-1}.$$
	This permutation induces a mapping from the two-dimensional index set
	$(l,m)\in\{0,1,\cdots,N_1-1\}\times\mathcal K_{N_2}$
	to the one-dimensional index set $\{0,1,\cdots,N_1N_2-1\}$, defined by
	\begin{equation*}
		r(l,m):=\iota^{-1}(m)N_1+l.
	\end{equation*}
	This mapping specifies a linear ordering of the two-dimensional indices. Hence, the Fourier coefficient vector \eqref{eq:U vector} can be identified with the standard vector $\widetilde {\bm U}^{(r)}=\big(\widetilde U_{q}^{(r)}\big)_{q=0}^{N_1N_2-1}$, satisfying
	$\widetilde U_{r(l,m)}^{(r)}=\widetilde U_{l,m}$. Under the same ordering, the preconditioner \eqref{eq:preconditioner P} has the matrix representation $\mathbf P^{(r)}=(P^{(r)}_{q_1,q_2})_{q_1,q_2=0}^{N_1N_2-1}$, where
	$$P^{(r)}_{r(l_1,m_1),r(l_2,m_2)}=P_{(l_1,m_1),(l_2,m_2)}.$$
	In view of the Fourier-mode shift relations at the boundaries \eqref{eq:ghost Fourier}, the matrix entry $P^{(r)}_{q_1,q_2}$ can be nonzero only in the following cases:
	\begin{enumerate}
		\item $q_1=q_2$;
		\item $|q_1-q_2|=1$, with $\Big\lfloor \dfrac{q_1d}{N_1N_2}\Big\rfloor=\Big\lfloor \dfrac{q_2d}{N_1N_2}\Big\rfloor$, where $\lfloor\cdot\rfloor$ denotes the floor function;
		\item $|q_1-q_2|=\dfrac{N_1N_2}{d}-1$, with $\Big\lfloor \dfrac{q_1d}{N_1N_2}\Big\rfloor=\Big\lfloor \dfrac{q_2d}{N_1N_2}\Big\rfloor$.
	\end{enumerate}
	The numbers of nonzero entries in these three categories are $N_1N_2$, $2N_1N_2-2d$, and $2d$, respectively. Only the third category lies outside the tridiagonal band.
	
	Since $\mathbf P^{(r)}$ is simply the matrix representation of $\mathbf P$ under the above permutation, it suffices to consider its efficient inversion. We construct a lower triangular factor $\mathbf L^{(r)}=(L_{q_1,q_2})_{q_1,q_2=0}^{N_1N_2-1}$ through the Cholesky factorization
	\begin{equation}\label{eq:cholesky}
		\mathbf P^{(r)}=\mathbf L^{(r)}\big(\mathbf L^{(r)}\big)^*.
	\end{equation}
	The entries of $\mathbf L^{(r)}$ can be computed using the recursions
	\begin{equation*}
		\begin{aligned}
			&L_{q_1,q_1}:=\bigg(P^{(r)}_{q_1,q_1}-\sum_{s=0}^{q_1-1}|L_{q_1,s}|^2\bigg)^{1/2},
			\qquad q_1=0,1,\ldots,N_1N_2-1,\\
			&L_{q_1,q_2}:=\frac{P^{(r)}_{q_1,q_2}-\sum_{s=0}^{q_2-1}L_{q_1,s}L_{q_2,s}^*}{L_{q_2,q_2}},
			\qquad 0\leq q_2<q_1\leq N_1N_2-1.
		\end{aligned}
	\end{equation*}
	The sparsity pattern of $\mathbf P^{(r)}$ implies that nonzero entries of the Cholesky factor $\mathbf L^{(r)}$ may occur only in the following three cases:
	\begin{enumerate}
		\item $q_1=q_2$;
		\item $q_1-q_2=1$, with $\Big\lfloor \dfrac{q_1d}{N_1N_2}\Big\rfloor=\Big\lfloor \dfrac{q_2d}{N_1N_2}\Big\rfloor$;
		\item $q_1\,\mathrm{MOD}\,\dfrac{N_1N_2}{d}=\dfrac{N_1N_2}{d}-1$, with $\Big\lfloor \dfrac{q_1d}{N_1N_2}\Big\rfloor=\Big\lfloor \dfrac{q_2d}{N_1N_2}\Big\rfloor$.
	\end{enumerate}
	The third case accounts for the additional fill-in generated by the Cholesky factorization. Since the total number of fill-in entries does not exceed $3N_1N_2$, the factor $\mathbf L^{(r)}$ remains sparse. Within the \gls{LOBPCG} iteration, the preconditioner is applied primarily to the residual $\widetilde{\bm R}$, requiring the computation of
	$$\widetilde{\bm Z}^{(r)}=(\mathbf P^{(r)})^{-1}\widetilde{\bm R}^{(r)},$$
	Given the Cholesky factorization \eqref{eq:cholesky}, this operation reduces to two triangular solves:
	$$\mathbf L^{(r)}\widetilde{\bm Y}^{(r)}=\widetilde{\bm R}^{(r)},\qquad \big(\mathbf L^{(r)}\big)^*\widetilde{\bm Z}^{(r)}=\widetilde{\bm Y}^{(r)}.$$
	Exploiting the sparsity of $\mathbf L^{(r)}$, each application of the preconditioner requires $O(N_1N_2)$ operations.

	\section{Error Analysis of the Fully Discrete Scheme}\setcounter{equation}{0}
	
	\subsection{Fourier-coefficient sampling and Sobolev spaces}
	
	Let $\mathcal{C}_B$ denote the space of continuous functions satisfying the magnetic Bloch boundary conditions, namely,
	$$\mathcal{C}_B:=\big\{u\in C(\mathbb{R}^2)\,:\,u(x+a,y)=\mathrm e^{\mathrm i\bar Bay}u(x,y),~u(x,y+b)=u(x,y)\big\}.$$
	To represent continuous functions by discrete Fourier coefficient vectors, we define the sampling operator \(\mathbf{\widetilde R}_{\bm h}:\mathcal C_B\to\mathbb C^{N_1N_2}\). This operator first samples the function values on the computational grid and then applies the \gls{DFT} in the $y$-direction, i.e.,
	\begin{equation*}
		\big(\mathbf{\widetilde{R}}_{\bm h}u\big)_{l,j}:=\frac1{N_2}\sum_{m=0}^{N_2-1}u(x_l,y_m)\mathrm e^{-\mathrm i\kappa_jy_m},
		\qquad j\in\mathcal K_{N_2}.
	\end{equation*}
	By Parseval's identity, we have
	\begin{equation}\label{eq:err parseval physical Fourier}
		\begin{aligned}
			h_1h_2\sum_{l=0}^{N_1-1}\sum_{m=0}^{N_2-1}
			u(x_l,y_m)
			\overline{w(x_l,y_m)}
			&=h_1b\sum_{l=0}^{N_1-1}
			\sum_{j\in\mathcal K_{N_2}}
			\big(\mathbf{\widetilde R}_{\bm h}u\big)_{l,j}
			\overline{\big(\mathbf{\widetilde R}_{\bm h}w\big)_{l,j}}\\
			&=\big\langle
			\mathbf{\widetilde R}_{\bm h}u,
			\mathbf{\widetilde R}_{\bm h}w
			\big\rangle_{\bm h}.
		\end{aligned}
	\end{equation}
	In particular, setting $w=u$ yields
	\[
	h_1h_2\sum_{l=0}^{N_1-1}\sum_{m=0}^{N_2-1}
	|u(x_l,y_m)|^2
	=
	\big\|\mathbf{\widetilde R}_{\bm h}u\big\|_{\bm h}^2.
	\]
	
	To characterize the regularity of the exact eigenfunctions, we introduce the anisotropic Sobolev norm
	\begin{equation*}
		\|u\|_{\mathcal{H}^{r,s}(\Omega)}^2
		:=b\sum_{j\in\mathbb Z}(1+\kappa_j^2)^{s}
		\|\widehat u_j\|_{\mathcal{H}^r(-a/2,a/2)}^2,
	\end{equation*}
	where $\|\cdot\|_{\mathcal{H}^r(-a/2,a/2)}$ denotes the standard Sobolev norm for functions of one variable.
	For functions with sufficiently high regularity, the corresponding Fourier coefficients also satisfy the following higher-order magnetic Bloch compatibility conditions:
	\begin{equation}\label{eq:err high order magnetic compatibility}
		\partial_x^\ell\widehat u_j\left(\frac a2\right)
		=\partial_x^\ell\widehat u_{j-n_B}\left(-\frac a2\right),
		\qquad \ell<r-\frac12.
	\end{equation}
	We can thus define the Sobolev space incorporating the magnetic Bloch boundary conditions by
	\begin{equation*}
		\mathcal{H}_B^{r,s}(\Omega):=\left\{u:\|u\|_{\mathcal{H}^{r,s}(\Omega)}<\infty,
		\text{and satisfying}\ \eqref{eq:err high order magnetic compatibility}\right\}.
	\end{equation*}

	\subsection{Second-order consistency of the covariant central difference in $x$}
	
	Denote the one-dimensional torus by
	\[
	\mathbb T_a:=\mathbb R/(a\mathbb Z).
	\]
	For a fixed $y=y_m$, let
	\begin{equation*}
		D_{x,m}:=-\mathrm i\frac{d}{dx}+k_1+\phi_y(x,y_m).
	\end{equation*}
	For the analysis below, define the Sobolev space
	\[
	\begin{aligned}
		\mathcal{H}_{B,m}^s(\mathbb T_a)
		:=
		\big\{
		w\in \mathcal{H}_{\rm loc}^s(\mathbb R):
		w^{(j)}(x+a)
		=
		\mathrm e^{\mathrm i\bar B ay_m} w^{(j)}(x),\quad
		j=0,1,\cdots,s-1
		\big\},
	\end{aligned}
	\]
	where $s\in\mathbb{N}_+$.
	
	Define the one-dimensional sampling operator and discrete norm by
	\[
	(\mathbf{R}_{h_1}w)_l:=w(x_l),
	\qquad
	\|\bm W\|_{h_1}^2
	:=
	h_1\sum_{l=0}^{N_1-1}|W_l|^2.
	\]
	We extend discrete vectors using the magnetic Bloch condition:
	\[
	W_{l+N_1}=\mathrm e^{\mathrm i\bar B ay_m}W_l,
	\qquad
	G_{l+N_1+1/2,m}=G_{l+1/2,m}.
	\]
	The one-dimensional covariant central difference operator is defined by
	\begin{equation}\label{eq:err one dimensional link operator}
		\left(\mathbf{L}_{x,h_1,m}\bm W\right)_l
		:=
		\frac{
			2W_l-G_{l+1/2,m}W_{l+1}
			-\overline{G_{l-1/2,m}}W_{l-1}
		}{h_1^2}.
	\end{equation}
	
	\begin{lemma}
		\label{lem:err x link consistency}
		Assume that
		\[
		\sup_{y\in\mathbb T_b}
		\|\phi_y(\cdot,y)\|_{\mathcal{W}^{3,\infty}(-a/2,a/2)}<\infty,
		\]
		and that $w\in\mathcal{H}_{B,m}^4(\mathbb T_a)$. Then
		\begin{equation}\label{eq:err x link consistency estimate}
			\left\|
			\mathbf{L}_{x,h_1,m}\mathbf{R}_{h_1}w
			-
			\mathbf{R}_{h_1}(D_{x,m}^2w)
			\right\|_{h_1}
			\leq
			Ch_1^2
			\|w\|_{\mathcal{H}^4(-a/2,a/2)},
		\end{equation}
		Furthermore, if $u\in\mathcal{H}_B^{4,1}(\Omega)$, then
		\begin{equation}\label{eq:err Ex global Sobolev bound}
			\big\|\mathbf L_{x,h_1}\mathbf{\widetilde{R}}_{\bm h}u-
			\mathbf{\widetilde{R}}_{\bm h}(D_x^2u)\big\|_{\bm h}
			\leq
			Ch_1^2\|u\|_{\mathcal{H}^{4,1}(\Omega)},
		\end{equation}
		The constant $C$ in these estimates is independent of $h_1$, $N_2$, and the norms of $w$ and $u$.
	\end{lemma}
	
	\begin{proof}
		First, consider the operator with exact link phases:
		$$(\mathbf{L}_{x,h_1,m}^{e}\bm W)_l
		:=
		\frac{
			2W_l-\mathrm e^{\mathrm i\vartheta_{l,m}(x_{l+1})}W_{l+1}
			-\mathrm e^{\mathrm i\vartheta_{l,m}(x_{l-1})}W_{l-1}
		}{h_1^2}.$$
		For any fixed grid index $l$, define
		\[
		v_l(x):=\mathrm e^{\mathrm i\vartheta_{l,m}(x)}w(x).
		\]
		Then, on the local stencil centered at $x_l$, we have
		\[
		D_{x,m}^2w(x)
		=-\mathrm e^{-\mathrm i\vartheta_{l,m}(x)}v_l''(x).
		\]
		Thus, the exact-link operator $\mathbf{L}_{x,h_1,m}^{e}$ reduces to the standard second-order central difference approximation of
		$-v_l''$; explicitly,
		\[
		\begin{aligned}
			D_{x,m}^2w(x_l)
			&\approx \mathrm e^{-\mathrm i\vartheta_{l,m}(x_l)}
			\frac{2v_l(x_l)-v_l(x_{l+1})-v_l(x_{l-1})}{h_1^2}\\
			&=\frac{
				2W_l-\mathrm e^{\mathrm i\vartheta_{l,m}(x_{l+1})}W_{l+1}
				-\mathrm e^{-\mathrm i\vartheta_{l-1,m}(x_l)}W_{l-1}
			}{h_1^2}.
		\end{aligned}
		\]
		Since
		$\sup_{y\in\mathbb T_b}\|\phi_y(\cdot,y)\|_{\mathcal{W}^{3,\infty}(-a/2,a/2)}<\infty$,
		the phase multiplier $\mathrm e^{\mathrm i\vartheta_{l,m}}$ is uniformly bounded on $\mathcal{H}^4(-a/2,a/2)$, namely,
		\[
		\|v_l\|_{\mathcal{H}^4(-a/2,a/2)}
		\leq C\|w\|_{\mathcal{H}^4(-a/2,a/2)}.
		\]
		Applying the standard discrete $L^2$ consistency estimate for the central difference scheme, we obtain
		\begin{equation}\label{eq:err exact link Sobolev consistency}
			\left\|
			\mathbf{L}_{x,h_1,m}^{e}\mathbf{R}_{h_1}w
			-
			\mathbf{R}_{h_1}(D_{x,m}^2w)
			\right\|_{h_1}
			\leq
			Ch_1^2
			\|w\|_{\mathcal{H}^4(-a/2,a/2)}.
		\end{equation}
		
		Let
		\[
		\delta_{l+1/2,m}
		:=
		\Theta_{l+1/2,m}
		-
		\vartheta_{l,m}(x_{l+1}).
		\]
		The trapezoidal-rule error and the assumption $\phi_y(\cdot,y)\in\mathcal{W}^{3,\infty}(\mathbb T_a)$ imply that
		\begin{equation}\label{eq:err trapezoid phase error}
			|\delta_{l+1/2,m}|
			\leq
			Ch_1^3\big\|\partial_x^2\phi_y(\cdot,y_m)\big\|_{\mathcal{L}^\infty(-a/2,a/2)},
		\end{equation}
		and
		\begin{equation}\label{eq:err phase error difference}
			|\delta_{l+1/2,m}-\delta_{l-1/2,m}|
			\leq
			Ch_1^4\big\|\partial_x^3\phi_y(\cdot,y_m)\big\|_{\mathcal{L}^\infty(-a/2,a/2)}.
		\end{equation}
		Using
		\[
		G_{l+1/2,m}
		=
		\mathrm e^{\mathrm i\vartheta_{l,m}(x_{l+1})}\mathrm e^{\mathrm i\delta_{l+1/2,m}},
		\]
		and the expansion $\mathrm e^{\mathrm i\delta}-1=\mathrm i\delta+O(\delta^2)$ in the definition \eqref{eq:err one dimensional link operator} of $\mathbf{L}_{x,h_1,m}$, we obtain
		\[
		\begin{aligned}
			&
			\left(
			\left(
			\mathbf{L}_{x,h_1,m}
			-
			\mathbf{L}_{x,h_1,m}^e
			\right)\mathbf{R}_{h_1}w
			\right)_{l}
			\\
			&\quad=
			-\frac{\mathrm i}{h_1^2}
			\left[
			\mathrm e^{\mathrm i\vartheta_{l,m}(x_{l+1})}
			\delta_{l+1/2,m}
			w(x_{l+1})
			-
			\mathrm e^{-\mathrm i\vartheta_{l-1,m}(x_{l})}
			\delta_{l-1/2,m}
			w(x_{l-1})
			\right]+R_l,
		\end{aligned}
		\]
		where the remainder terms $R_l$, $l=0,1,\ldots,N_1-1$, satisfy
		\[
		\big\|(R_l)_{l=0}^{N_1-1}\big\|_{h_1}
		\leq
		Ch_1^3
		\|w\|_{\mathcal{H}^1(-a/2,a/2)}.
		\]
		By \eqref{eq:err trapezoid phase error} and \eqref{eq:err phase error difference}, the leading term can be bounded as follows:
		\[
		\begin{aligned}
			&
			\frac{1}{h_1^2}
			\left|
			\mathrm e^{\mathrm i\vartheta_{l,m}(x_{l+1})}
			\delta_{l+1/2,m}
			w(x_{l+1})
			-
			\mathrm e^{-\mathrm i\vartheta_{l-1,m}(x_{l})}
			\delta_{l-1/2,m}
			w(x_{l-1})
			\right|
			\\
			&\quad\leq
			\frac{|\delta_{l+1/2,m}
				-\delta_{l-1/2,m}|}{h_1^2}
			|w(x_{l+1})|\\
			&\qquad\quad+
			\frac{|\delta_{l-1/2,m}|}{h_1^2}
			\left|
			\mathrm e^{\mathrm i\vartheta_{l,m}(x_{l+1})}
			w(x_{l+1})
			-
			\mathrm e^{\mathrm i\vartheta_{l,m}(x_{l-1})}
			w(x_{l-1})
			\right|\\
			&\quad\leq Ch_1^2\|w\|_{\mathcal{W}^{1,\infty}(-a/2,a/2)}.
		\end{aligned}
		\]
		The one-dimensional Sobolev embedding
		$\|w\|_{\mathcal{W}^{1,\infty}(-a/2,a/2)}
		\leq C\|w\|_{\mathcal{H}^2(-a/2,a/2)}$ then yields
		\begin{equation}\label{eq:err trapezoid exact L2}
			\left\|
			\left(
			\mathbf{L}_{x,h_1,m}
			-
			\mathbf{L}_{x,h_1,m}^e
			\right)\mathbf{R}_{h_1}w
			\right\|_{h_1}
			\leq
			Ch_1^2
			\|w\|_{\mathcal{H}^4(-a/2,a/2)}.
		\end{equation}
		Combining \eqref{eq:err exact link Sobolev consistency} with
		\eqref{eq:err trapezoid exact L2} gives
		\eqref{eq:err x link consistency estimate}.
		
		For a two-dimensional physical-space grid function $\bm U=(U_{l,m})$, define
		\begin{equation*}
			\mathbf{L}_{x,h_1}^{P}\bm U
			:=
			\left(\mathbf{L}_{x,h_1,m}(U_{l,m})_{l=0}^{N_1-1}\right)_{m=0}^{N_2-1},
		\end{equation*}
		where the operator $\mathbf{L}_{x,h_1}^{P}$ is related to $\mathbf L_{x,h_1}$ by
		\[
		\mathbf L_{x,h_1}
		=
		\mathbf F_y\mathbf{L}_{x,h_1}^{P}\mathbf F_y^{-1}.
		\]
		By the discrete Parseval identity \eqref{eq:parseval},
		\[
		\begin{aligned}
			\big\|\mathbf L_{x,h_1}\mathbf{\widetilde{R}}_{\bm h}u-
			\mathbf{\widetilde{R}}_{\bm h}(D_x^2u)\big\|_{\bm h}^2
			&=
			h_2\sum_{m=0}^{N_2-1}
			\left\|
			\mathbf{L}_{x,h_1,m}\mathbf{R}_{h_1}u(\cdot,y_m)
			-
			\mathbf{R}_{h_1}D_{x,m}^2u(\cdot,y_m)
			\right\|_{h_1}^2\\
			&\leq
			Ch_1^4h_2\sum_{m=0}^{N_2-1}
			\|u(\cdot,y_m)\|_{\mathcal{H}^4(-a/2,a/2)}^2.
		\end{aligned}
		\]
		Finally, using the vector-valued sampling stability estimate
		\[
		h_2\sum_{m=0}^{N_2-1}
		\|u(\cdot,y_m)\|_{\mathcal{H}^4(-a/2,a/2)}^2
		\leq
		C\|u\|_{\mathcal{H}^{4,1}(\Omega)}^2,
		\]
		we obtain \eqref{eq:err Ex global Sobolev bound}.
	\end{proof}

	\subsection{Variable-coefficient collocation errors in Fourier pseudospectral discretization}
	This subsection analyzes the consistency error of the Fourier pseudospectral discretization in the $y$-direction, including Fourier interpolation errors and collocation errors arising from multiplication by the variable coefficient $\phi_x$. The analysis relies primarily on Fourier interpolation estimates and multiplication stability in periodic Sobolev spaces; see, for example, \cite{LiXiao2021}.
	
	Let $\mathcal I_{h_2}$ be the Fourier interpolation operator defined by
	$$\mathcal I_{h_2}u(x,y):=\sum_{j\in\mathcal{K}_{N_2}}\widetilde{u}_{j}(x)\mathrm e^{\mathrm i\kappa_jy},$$
	where $\big(\widetilde{u}_{j}(x)\big)_{j\in\mathcal{K}_{N_2}}=\mathbf{F}_y\big(u(x,y_m)\big)_{m=0}^{N_2-1}$.
	For a fixed $x=x_l$, the covariant derivative in the $y$-direction takes the form
	\begin{equation*}
		D_{y,l}w(y):=-\mathrm i\partial_yw(y)+k_2w(y)-\bar Bx_lw(y)-\phi_x(x_l,y)w(y).
	\end{equation*}
	The discrete counterpart of $D_{y,l}$ is given by
	\begin{equation*}
		\mathcal D_{y,h_2,l}w(y)
		:=-\mathrm i\partial_y w(y)
		+\mathcal I_{h_2}\big((k_2-\bar Bx_l-\phi_x(x_l,y))w(y)\big),
	\end{equation*}
	where
	$w\in {\rm span}\{\mathrm e^{\mathrm i\kappa_jy}:j\in\mathcal K_{N_2}\}$.
	Consequently, $\mathbf L_{y,h_2}$ defined in \eqref{eq:Ly square} is the Fourier coefficient-vector representation of
	$\mathcal D_{y,h_2,l}^2$.
	
	\begin{lemma}\label{lem:err y variable coefficient concrete}
		Let $s>1/2$, $u\in\mathcal{H}_B^{4,s+2}(\Omega)$, and suppose that
		$\Phi_{s+2}:=\sup_{x\in (-\frac{a}{2},\frac{a}{2})} \|k_2-\bar Bx-\phi_x(x,\cdot)\|_{\mathcal{H}^{s+2}(\mathbb{T}_b)}<\infty$. Then
		\begin{equation}\label{eq:err y variable coefficient concrete}
			\big\|\mathbf L_{y,h_2}\mathbf{\widetilde{R}}_{\bm h}u-\mathbf{\widetilde{R}}_{\bm h}(D_y^2u)\big\|_{\bm h}
			\leq C N_2^{-s}\|u\|_{\mathcal{H}^{1,s+2}(\Omega)},
		\end{equation}
		The constant $C$ is independent of $h_1$, $N_2$, and the norm of $u$.
	\end{lemma}
	
	\begin{proof}
		Fix $x_l$ and set $u_l(y)=u(x_l,y)$ and $\varphi_l(y)=k_2-\bar Bx_l-\phi_x(x_l,y)$. By definition, the $l$th block of $\mathbf{\widetilde{R}}_{\bm h}u$ corresponds to the trigonometric polynomial $\mathcal I_{h_2}u_l$, while the continuous operator satisfies
		\begin{equation*}
			D_{y,l}^2u_l
			=-\partial_y^2u_l-\mathrm i\partial_y(\varphi_lu_l)-\mathrm i\varphi_l\partial_yu_l+\varphi_l^2u_l.
		\end{equation*}
		The collocation-discretized operator, in function form, can be written as
		\begin{equation*}
			\begin{aligned}
				\mathcal D_{y,h_2,l}^2(\mathcal I_{h_2}u_l)
				=&-\partial_y^2 \mathcal I_{h_2}u_l
				-\mathrm i\partial_y \mathcal I_{h_2}(\varphi_l\mathcal I_{h_2}u_l) -\mathrm i \mathcal I_{h_2}(\varphi_l\partial_y\mathcal I_{h_2}u_l)\\
				&+\mathcal I_{h_2}\big(\varphi_l\mathcal I_{h_2}(\varphi_l\mathcal I_{h_2}u_l)\big).
			\end{aligned}
		\end{equation*}
		Therefore,
		\begin{equation}\label{eq:err y residual four terms}
			\begin{aligned}
				\mathcal D_{y,h_2,l}^2(\mathcal I_{h_2}u_l)-D_{y,l}^2u_l
				=&-\partial_y^2(\mathcal I_{h_2}u_l-u_l)-\mathrm i\partial_y\big(\mathcal I_{h_2}(\varphi_l\mathcal I_{h_2}u_l)-\varphi_lu_l\big)                    \\
				&-i\big(\mathcal I_{h_2}(\varphi_l\partial_y\mathcal I_{h_2}u_l)-\varphi_l\partial_yu_l\big)          \\
				&+\big(\mathcal I_{h_2}(\varphi_l\mathcal I_{h_2}(\varphi_l\mathcal I_{h_2}u_l))-\varphi_l^2u_l\big).
			\end{aligned}
		\end{equation}
		We estimate the four terms separately.
		
		The Fourier interpolation error estimate states that, for $0\leq\gamma\leq\beta$ with $\beta>1/2$,
		\begin{equation}\label{eq:err interpolation estimate y}
			\|w-\mathcal I_{h_2}w\|_{\mathcal{H}^\gamma(\mathbb{T}_b)}
			\leq C N_2^{\gamma-\beta}\|w\|_{\mathcal{H}^\beta(\mathbb{T}_b)},
		\end{equation}
		see Eq.~(3.3) in \cite{LiXiao2021} for details. Taking $\gamma=2$ and $\beta=s+2$, we obtain
		\begin{equation}\label{eq:err term y second derivative}
			\big\|\partial_y^2(\mathcal I_{h_2}u_l-u_l)\big\|_{L^2(\mathbb{T}_b)}
			\leq C N_2^{-s}\|u_l\|_{\mathcal{H}^{s+2}(\mathbb{T}_b)}.
		\end{equation}
		
		For $\rho>1/2$, the one-dimensional periodic Sobolev space $\mathcal{H}^\rho(\mathbb{T}_b)$ is a Banach algebra and satisfies the product estimate
		\begin{equation}\label{eq:err Sobolev product y}
			\|fg\|_{\mathcal{H}^\rho(\mathbb{T}_b)}
			\leq C\|f\|_{\mathcal{H}^\rho(\mathbb{T}_b)}
			\|g\|_{\mathcal{H}^\rho(\mathbb{T}_b)},
			\qquad \rho>1/2,
		\end{equation}
		Combining this with the coefficient regularity assumption $\Phi_{s+2}<\infty$ gives
		\begin{equation*}
			\begin{aligned}
				&\|\varphi_lu_l\|_{\mathcal{H}^{s+1}(\mathbb{T}_b)}
				+\|\varphi_l\partial_yu_l\|_{\mathcal{H}^{s}(\mathbb{T}_b)}
				+\|\varphi_l^2u_l\|_{\mathcal{H}^{s}(\mathbb{T}_b)}\\
				&\quad\leq C(1+\Phi_{s+2}^2)\|u_l\|_{\mathcal{H}^{s+2}(\mathbb{T}_b)}.
			\end{aligned}
		\end{equation*}
		Moreover, $\mathcal I_{h_2}$ is stable in $\mathcal{H}^\rho(\mathbb{T}_b)$, and \eqref{eq:err interpolation estimate y} implies
		\begin{equation*}
			\|u_l-\mathcal I_{h_2}u_l\|_{\mathcal{H}^{1}(\mathbb{T}_b)}
			\leq C N_2^{-(1+s)}\|u_l\|_{\mathcal{H}^{s+2}(\mathbb{T}_b)}.
		\end{equation*}
		Hence,
		\begin{equation}\label{eq:err term derivative au}
			\begin{aligned}
				&\|\partial_y(\mathcal I_{h_2}(\varphi_l\mathcal I_{h_2}u_l)-\varphi_lu_l)\|_{L^2(\mathbb{T}_b)}       \\
				&\quad\leq
				\|\varphi_l\mathcal I_{h_2}u_l-\mathcal I_{h_2}(\varphi_l\mathcal I_{h_2}u_l)\|_{\mathcal{H}^1(\mathbb{T}_b)}
				+\|\varphi_l(\mathcal I_{h_2}u_l-u_l)\|_{\mathcal{H}^1(\mathbb{T}_b)}                         \\
				&\quad\leq C N_2^{-s}(1+\Phi_{s+2})\|u_l\|_{\mathcal{H}^{s+2}(\mathbb{T}_b)}.
			\end{aligned}
		\end{equation}
		Similarly, splitting the error into the Fourier interpolation residual and the interpolation error due to
		$\partial_y(\mathcal I_{h_2}u_l-u_l)$, we obtain
		\begin{equation}\label{eq:err term a uy}
			\begin{aligned}
				&\|\mathcal I_{h_2}(\varphi_l\partial_y\mathcal I_{h_2}u_l)-\varphi_l\partial_yu_l\|_{L^2(\mathbb{T}_b)}  \\
				&\quad\leq
				\|\varphi_l\partial_y\mathcal I_{h_2}u_l-\mathcal I_{h_2}(\varphi_l\partial_y\mathcal I_{h_2}u_l)\|_{L^2(\mathbb{T}_b)}
				+\|\varphi_l\partial_y(\mathcal I_{h_2}u_l-u_l)\|_{L^2(\mathbb{T}_b)}                 \\
				&\quad\leq C N_2^{-s}(1+\Phi_{s+2})\|u_l\|_{\mathcal{H}^{s+2}(\mathbb{T}_b)}.
			\end{aligned}
		\end{equation}
		For the double multiplication term, we first write
		\begin{equation*}
			\begin{aligned}
				&\mathcal I_{h_2}(\varphi_l\mathcal I_{h_2}(\varphi_l\mathcal I_{h_2}u_l))-\varphi_l^2u_l\\
				&\quad=\mathcal I_{h_2}\big(\varphi_l\mathcal I_{h_2}(\varphi_l\mathcal I_{h_2}u_l)\big)-\varphi_l\mathcal I_{h_2}(\varphi_l\mathcal I_{h_2}u_l)
				+\varphi_l\left[\mathcal I_{h_2}(\varphi_l\mathcal I_{h_2}u_l)-\varphi_lu_l\right].
			\end{aligned}
		\end{equation*}
		Applying \eqref{eq:err interpolation estimate y} and \eqref{eq:err Sobolev product y} once more yields
		\begin{equation}\label{eq:err term a2u}
			\|\mathcal I_{h_2}(\varphi_l\mathcal I_{h_2}(\varphi_l\mathcal I_{h_2}u_l))-\varphi_l^2u_l\|_{L^2(\mathbb{T}_b)}
			\leq C N_2^{-s}(1+\Phi_{s+2}^2)\|u_l\|_{\mathcal{H}^{s+2}(\mathbb{T}_b)}.
		\end{equation}
		Substituting \eqref{eq:err term y second derivative}, \eqref{eq:err term derivative au},
		\eqref{eq:err term a uy}, and \eqref{eq:err term a2u} into \eqref{eq:err y residual four terms}, we obtain, for each $l$,
		\begin{equation*}
			\begin{aligned}
				&\|\mathcal D_{y,h_2,l}^2(\mathcal I_{h_2}u_l)-\mathcal I_{h_2}(D_{y,l}^2u_l)\|_{L^2(\mathbb{T}_b)}\\
				&\quad \leq\|\mathcal D_{y,h_2,l}^2(\mathcal I_{h_2}u_l)-D_{y,l}^2u_l\|_{L^2(\mathbb{T}_b)}+\|D_{y,l}^2u_l-\mathcal I_{h_2}(D_{y,l}^2u_l)\|_{L^2(\mathbb{T}_b)}\\
				&\quad\leq C N_2^{-s}\|u_l\|_{\mathcal{H}^{s+2}(\mathbb{T}_b)}.
			\end{aligned}
		\end{equation*}
		Finally, summing over the grid points $x_l$ and applying the sampling stability estimate in the $x$-direction, we have
		\begin{equation*}
			h_1\sum_{l=0}^{N_1-1}\|u(x_l,\cdot)\|_{\mathcal{H}^{s+2}(\mathbb{T}_b)}^2
			\leq C\|u\|_{\mathcal{H}^{1,s+2}(\Omega)}^2.
		\end{equation*}
		The discrete Parseval identity \eqref{eq:parseval} then yields
		\eqref{eq:err y variable coefficient concrete}, completing the proof.
	\end{proof}

	\subsection{Residual estimates and errors for discrete eigenvalue clusters}
	Let $\lambda^{(\bm{k})}$ be an isolated eigenvalue of the continuous Hamiltonian $H^{(\bm{k})}$,
	with associated eigenfunction $u$, so that
	\begin{equation}\label{eq:err continuous eigen equation}
		H^{(\bm{k})}u=\lambda^{(\bm{k})} u.
	\end{equation}
	
	We estimate the consistency residual obtained by sampling an eigenfunction on the grid, applying the \gls{DFT} in the $y$-direction, and inserting the resulting coefficient vector into the discrete eigenvalue equation.
	
	\begin{lemma}
		\label{lem:err residual consistency link}
		Let $s>1/2$ and $u\in\mathcal{H}_B^{4,s+2}(\Omega)$.
		Assume that $\phi_y$ and $\phi_x$ satisfy the respective regularity conditions of Lemmas \ref{lem:err x link consistency} and \ref{lem:err y variable coefficient concrete}.
		If $B$, $V$, $\phi_x$, and $\phi_y$ are all sampled exactly from their continuous counterparts in the discrete scheme,
		then
		\begin{equation}\label{eq:err residual consistency one}
			\begin{aligned}
				\big\|
				\big(\mathbf H_{\bm h}^{(\bm{k})}-\lambda^{(\bm{k})}\mathbf I\big)
				\mathbf{\widetilde R}_{\bm h}u
				\big\|_{\bm h}
				\leq
				C\left(h_1^2+N_2^{-s}\right)
				\|u\|_{\mathcal{H}^{4,s+2}(\Omega)}.
			\end{aligned}
		\end{equation}
		where the constant $C$ is independent of $h_1$, $N_2$, and the norm of $u$.
	\end{lemma}
	
	\begin{proof}
		The continuous eigenvalue equation \eqref{eq:err continuous eigen equation} implies
		\begin{equation*}
			\big(\mathbf H_{\bm h}^{(\bm{k})}-\lambda^{(\bm{k})}\mathbf I\big)
			\mathbf{\widetilde R}_{\bm h}u=
			\mathbf H_{\bm h}^{(\bm{k})}
			\mathbf{\widetilde R}_{\bm h}u
			-
			\mathbf{\widetilde R}_{\bm h}
			\left(H^{(\bm{k})}u\right).
		\end{equation*}
		By the definition \eqref{eq:Hamilton discrete} of the discrete Hamiltonian, the right-hand side can be decomposed as
		\begin{equation}\label{eq:err residual decomposition}
			\mathbf H_{\bm h}^{(\bm{k})}
			\mathbf{\widetilde R}_{\bm h}u
			-
			\mathbf{\widetilde R}_{\bm h}
			\left(H^{(\bm{k})}u\right)=
			\frac1{2}
			\big(
			\widetilde{\bm E}_{x,\bm h}(u)
			+
			\widetilde{\bm E}_{y,\bm h}(u)
			\big)
			+
			\widetilde{\bm E}_{V,\bm h}(u),
		\end{equation}
		where
		\[
		\begin{aligned}
			\widetilde{\bm E}_{x,\bm h}(u)
			&:=
			\mathbf L_{x,h_1}\mathbf{\widetilde R}_{\bm h}u
			-
			\mathbf{\widetilde R}_{\bm h}(D_x^2u),\\
			\widetilde{\bm E}_{y,\bm h}(u)
			&:=
			\mathbf L_{y,h_2}\mathbf{\widetilde R}_{\bm h}u
			-
			\mathbf{\widetilde R}_{\bm h}(D_y^2u),\\
			\widetilde{\bm E}_{V,\bm h}(u)
			&:=
			\mathbf{\widetilde V}_{\bm h}\mathbf{\widetilde R}_{\bm h}u
			-
			\mathbf{\widetilde R}_{\bm h}(Vu).
		\end{aligned}
		\]
		By Lemma \ref{lem:err x link consistency} and
		$\mathcal{H}_B^{4,s+2}(\Omega)\subset\mathcal{H}_B^{4,1}(\Omega)$,
		the error in the $x$-direction satisfies
		\begin{equation}\label{eq:err Ex bound}
			\big\|\widetilde{\bm E}_{x,\bm h}(u)\big\|_{\bm h}
			\leq
			Ch_1^2\|u\|_{\mathcal{H}^{4,s+2}(\Omega)}.
		\end{equation}
		By Lemma \ref{lem:err y variable coefficient concrete} and the embedding
		$\mathcal{H}_B^{4,s+2}(\Omega)\subset\mathcal{H}_B^{1,s+2}(\Omega)$, the error in the $y$-direction satisfies
		\begin{equation*}
			\big\|\widetilde{\bm E}_{y,\bm h}(u)\big\|_{\bm h}
			\leq
			CN_2^{-s}\|u\|_{\mathcal{H}^{4,s+2}(\Omega)}.
		\end{equation*}
		When the potential is sampled exactly at the grid points, both the discrete potential term and the sampled continuous product are obtained by forming $V(x_l,y_m)u(x_l,y_m)$ at the grid points and then applying the \gls{DFT} in the $y$-direction. The two expressions therefore coincide, namely,
		\begin{equation*}
			\mathbf{\widetilde V}_{\bm h}\mathbf{\widetilde R}_{\bm h}u
			=
			\mathbf{\widetilde R}_{\bm h}(Vu),
		\end{equation*}
		Consequently,
		\begin{equation}\label{eq:err EV bound}
			\big\|\widetilde{\bm E}_{V,\bm h}(u)\big\|_{\bm h}=0.
		\end{equation}
		Substituting \eqref{eq:err Ex bound}--\eqref{eq:err EV bound}
		into \eqref{eq:err residual decomposition} yields
		\eqref{eq:err residual consistency one}, as required.
	\end{proof}
	
	Let $\lambda^{(\bm{k})}$ be an isolated eigenvalue of multiplicity $m$. Its associated continuous eigenspace is
	\begin{equation*}
		\mathcal M(\lambda^{(\bm{k})})
		:=
		\ker\left(H^{(\bm{k})}-\lambda^{(\bm{k})} I\right),
		\qquad
		\dim\mathcal M(\lambda^{(\bm{k})})=m.
	\end{equation*}
	Since the eigenvalue is isolated, its distance from the remainder of the spectrum is positive. Denote this spectral gap by
	\begin{equation}\label{eq:err spectral gap}
		\delta_{\lambda^{(\bm{k})}}
		:=
		\operatorname{dist}\left(
		\lambda^{(\bm{k})},\,
		\sigma(H^{(\bm{k})})\setminus\{\lambda^{(\bm{k})}\}
		\right)>0.
	\end{equation}
	Choose an $L^2(\Omega)$-orthonormal basis of $\mathcal M(\lambda^{(\bm{k})})$,
	$u_1,\ldots,u_m$. %and define the sampled vectors
	%\begin{equation}\label{eq:err sampled vectors}
	%    \widetilde{\bm z}_{i,\bm h}
	%    :=
	%    \mathbf{\widetilde R}_{\bm h}u_i,
	%    \qquad i=1,\ldots,m.
	%\end{equation}
	The corresponding sampled Gram matrix is defined by
	\begin{equation*}
		\mathbf G_{\bm h}
		:=
		\left(
		\big\langle
		\mathbf{\widetilde R}_{\bm h}u_i,
		\mathbf{\widetilde R}_{\bm h}u_j
		\big\rangle_{\bm h}
		\right)_{i,j=1}^m.
	\end{equation*}
	For a matrix $\mathbf Z$ mapping $\mathbb C^m$ into the discrete Fourier coefficient space,
	define the associated induced norm by
	\begin{equation*}
		\|\mathbf Z\|_{\bm h,2}
		:=
		\sup_{\|\bm c\|_{2}=1}
		\|\mathbf Z\bm c\|_{\bm h}.
	\end{equation*}
	
	We estimate the deviation of the sampled Gram matrix from the identity, which will be used to construct a discrete orthonormal basis.
	
	\begin{lemma}
		\label{lem:err Gram consistency}
		Let $u_1,\ldots,u_m\in\mathcal{H}_B^{4,s+2}(\Omega)$ be an $L^2(\Omega)$-orthonormal basis of $\mathcal M(\lambda^{(\bm{k})})$, where $s>1/2$. Then
		\begin{equation*}
			\|\mathbf G_{\bm h}-\mathbf I_m\|_{2}
			\leq
			C\left(h_1^2+N_2^{-s}\right).
		\end{equation*}
		In particular, if $h_1$ is sufficiently small and $N_2$ is sufficiently large,
		$\mathbf G_{\bm h}$ is Hermitian positive definite, and
		$\|\mathbf G_{\bm h}^{-1/2}\|_2$ is uniformly bounded.
	\end{lemma}
	
	\begin{proof}
		By Parseval's identity \eqref{eq:err parseval physical Fourier},
		$$\big\langle
		\mathbf{\widetilde R}_{\bm h}u_i,
		\mathbf{\widetilde R}_{\bm h}u_j
		\big\rangle_{\bm h}=h_1h_2\sum_{l=0}^{N_1-1}\sum_{m=0}^{N_2-1}
		u_i(x_l,y_m)
		\overline{u_j(x_l,y_m)}.$$
		Since $u_i(x,y)$ and $u_j(x,y)$ satisfy the same magnetic Bloch boundary conditions,
		the product $u_i(x,y)\overline{u_j(x,y)}$ is ordinarily periodic in both the $x$- and $y$-directions.
		Combining the second-order periodic trapezoidal quadrature error in $x$ with the Fourier interpolation error in $y$ gives
		\[
		\big\langle
		\mathbf{\widetilde R}_{\bm h}u_i,
		\mathbf{\widetilde R}_{\bm h}u_j
		\big\rangle_{\bm h}
		=
		\langle u_i,u_j\rangle
		+
		O\left(h_1^2+N_2^{-s}\right).
		\]
		The conclusion follows from the orthonormality of $u_1,\ldots,u_m$.
	\end{proof}
	
	Define the discrete orthonormal basis matrix by
	\begin{equation}\label{eq:err Qh definition}
		\mathbf Q_{\bm h}
		:=
		\big[
		\mathbf{\widetilde R}_{\bm h}u_1,
		\ldots,
		\mathbf{\widetilde R}_{\bm h}u_m
		\big]
		\mathbf G_{\bm h}^{-1/2}.
	\end{equation}
	Let $\widetilde{\bm q}_{i,\bm h}$ denote the $i$th column of $\mathbf Q_{\bm h}$.
	Then
	\begin{equation*}
		\left(
		\big\langle
		\widetilde{\bm q}_{i,\bm h},
		\widetilde{\bm q}_{j,\bm h}
		\big\rangle_{\bm h}
		\right)_{i,j=1}^m
		=
		\mathbf I_m.
	\end{equation*}
	
	\begin{lemma}
		\label{lem:err subspace residual}
		Under the assumptions of Lemma \ref{lem:err residual consistency link}, we have
		\begin{equation}\label{eq:err subspace residual}
			\begin{aligned}
				\big\|
				\big(\mathbf H_{\bm h}^{(\bm{k})}-\lambda^{(\bm{k})}\mathbf I_{N_1N_2}\big)
				\mathbf Q_{\bm h}
				\big\|_{\bm h,2}
				\leq
				C\left(h_1^2+N_2^{-s}\right)
				\max_{1\leq i\leq m}
				\|u_i\|_{\mathcal{H}^{4,s+2}(\Omega)}.
			\end{aligned}
		\end{equation}
	\end{lemma}
	
	\begin{proof}
		By \eqref{eq:err Qh definition},
		\[
		\begin{aligned}
			&\big(\mathbf H_{\bm h}^{(\bm{k})}-\lambda^{(\bm{k})}\mathbf I_{N_1N_2}\big)
			\mathbf Q_{\bm h}\\
			&\quad=
			\big[
			\big(\mathbf H_{\bm h}^{(\bm{k})}-\lambda^{(\bm{k})}\mathbf I_{N_1N_2}\big)
			\mathbf{\widetilde R}_{\bm h}u_1,
			\ldots,
			\big(\mathbf H_{\bm h}^{(\bm{k})}-\lambda^{(\bm{k})}\mathbf I_{N_1N_2}\big)
			\mathbf{\widetilde R}_{\bm h}u_m
			\big]
			\mathbf G_{\bm h}^{-1/2}.
		\end{aligned}
		\]
		Applying Lemma \ref{lem:err residual consistency link} to each column,
		and using the fact that $m$ is fixed together with
		$\|\mathbf G_{\bm h}^{-1/2}\|_{2}\leq C$,
		we obtain \eqref{eq:err subspace residual}.
	\end{proof}
	
	We next use spectral projection estimates for finite-dimensional Hermitian matrices to relate the residual bound to eigenvalue and eigenspace errors; see \cite{BabuskaOsborn1991} for general spectral approximation theory. We assume that the target discrete eigenvalue cluster has the same total multiplicity as the isolated continuous eigenvalue.
	
	\begin{assumption}
		\label{ass:err discrete cluster dimension}
		Let $\lambda^{(\bm{k})}$ be an isolated eigenvalue of multiplicity $m$ of the continuous Hamiltonian
		$H^{(\bm{k})}$, and define
		\[
		I_{\delta_{\lambda^{(\bm{k})}}}
		:=
		\left(
		\lambda^{(\bm{k})}-\frac{\delta_{\lambda^{(\bm{k})}}}{2},
		\lambda^{(\bm{k})}+\frac{\delta_{\lambda^{(\bm{k})}}}{2}
		\right),
		\]
		where $\delta_{\lambda^{(\bm{k})}}$ is the spectral gap defined in
		\eqref{eq:err spectral gap}.
		Assume that there exist $h_{1,0},h_{2,0}>0$ such that, whenever
		
		\[
		0<h_1\leq h_{1,0},
		\qquad
		0<h_2\leq h_{2,0}
		\]
		the discrete Hermitian matrix
		$\mathbf H_{\bm h}^{(\bm{k})}$
		has exactly $m$ eigenvalues in the interval $I_{\delta_{\lambda^{(\bm{k})}}}$,
		counted with multiplicity.
	\end{assumption}
	
	By the dimension-preservation assumption for the discrete eigenvalue cluster, $\mathbf H_{\bm h}^{(\bm{k})}$ has exactly $m$ eigenpairs $(\lambda_{j,\bm h}^{(\bm{k})},\widetilde{\bm U}_{j,\bm h})$, $j=1,2,\ldots,m$, in $I_{\delta_{\lambda^{(\bm{k})}}}$, where $$\lambda^{(\bm{k})}-\delta_{\lambda^{(\bm{k})}}/2< \lambda_{1,\bm h}^{(\bm{k})}\leq\lambda_{2,\bm h}^{(\bm{k})}\leq\cdots\leq\lambda_{m,\bm h}^{(\bm{k})}<\lambda^{(\bm{k})}+\delta_{\lambda^{(\bm{k})}}/2,$$
	and the eigenvectors $\widetilde{\bm U}_{1,\bm h},\widetilde{\bm U}_{2,\bm h},\ldots,\widetilde{\bm U}_{m,\bm h}$ are orthonormal. Set
	$$\mathcal M_{\bm h}(\lambda^{(\bm{k})}):=\mathrm{span}\big\{\widetilde{\bm U}_{1,\bm h},\widetilde{\bm U}_{2,\bm h},\cdots,\widetilde{\bm U}_{m,\bm h}\big\}.$$
	The discrete spectral projector of $\mathbf H_{\bm h}^{(\bm{k})}$
	associated with $I_{\delta_{\lambda^{(\bm{k})}}}$, denoted by $\mathbf P_{\bm h}\,:\,\mathbb{C}^{N_1N_2}\rightarrow\mathcal M_{\bm h}(\lambda^{(\bm{k})})$, is defined by
	$$\mathbf P_{\bm h}\widetilde{\bm W}=\sum_{j=1}^m\big<\widetilde{\bm W},\widetilde{\bm U}_{j,\bm h}\big>_{\bm h}\widetilde{\bm U}_{j,\bm h}.$$
	It follows that
	$$\operatorname{Ran}\mathbf P_{\bm h}=\mathcal M_{\bm h}(\lambda^{(\bm{k})}),\qquad\operatorname{Ran}(\mathbf I-\mathbf P_{\bm h})=\mathcal M_{\bm h}(\lambda^{(\bm{k})})^{\bot}.$$
	
	\begin{theorem}
		\label{thm:err cluster}
		Let $\lambda^{(\bm{k})}$ be an isolated eigenvalue of multiplicity $m$ of $H^{(\bm{k})}$ satisfying the spectral gap condition
		\eqref{eq:err spectral gap}.
		If the discrete scheme $\mathbf H_{\bm h}^{(\bm{k})}$ satisfies Assumption \ref{ass:err discrete cluster dimension}, then
		\begin{equation}\label{eq:err eigenvalue cluster bound}
			\max_{1\leq i\leq m}
			|\lambda_{i,\bm h}^{(\bm{k})}-\lambda^{(\bm{k})}|
			\leq
			C\left(h_1^2+N_2^{-s}\right),
		\end{equation}
		and the sampled and discrete eigenspaces satisfy
		\begin{equation}\label{eq:err subspace gap bound}
			\left\|
			\left(\mathbf I-\mathbf P_{\bm h}\right)
			\mathbf Q_{\bm h}
			\right\|_{\bm h,2}
			\leq
			\frac{C}{\delta_{\lambda^{(\bm{k})}}}
			\left(h_1^2+N_2^{-s}\right).
		\end{equation}
		The constant $C$ may depend on the $\mathcal{H}^{4,s+2}(\Omega)$ norms of the basis functions of the target continuous eigenspace,
		but is independent of $h_1$ and $N_2$.
	\end{theorem}
	
	\begin{proof}
		Lemma \ref{lem:err subspace residual} yields
		\begin{equation}\label{eq:err Ah residual proof}
			\big\|
			\big(\mathbf H_{\bm h}^{(\bm{k})}-\lambda^{(\bm{k})}\mathbf I\big)
			\mathbf Q_{\bm h}
			\big\|_{\bm h,2}
			\leq
			C\left(h_1^2+N_2^{-s}\right).
		\end{equation}
		Since $\mathbf H_{\bm h}^{(\bm{k})}$ is Hermitian,
		its restriction to $\operatorname{Ran}(\mathbf I-\mathbf P_{\bm h})$ satisfies
		\begin{equation*}
			\big\|
			\big(\mathbf H_{\bm h}^{(\bm{k})}-\lambda^{(\bm{k})}\mathbf I\big)
			\widetilde{\bm W}
			\big\|_{\bm h}
			\geq
			\frac{\delta_{\lambda^{(\bm{k})}}}{2}
			\big\|\widetilde{\bm W}\big\|_{\bm h},
			\qquad
			\widetilde{\bm W}
			\in
			\operatorname{Ran}(\mathbf I-\mathbf P_{\bm h}).
		\end{equation*}
		Consequently,
		\[
		\begin{aligned}
			\frac{\delta_{\lambda^{(\bm{k})}}}{2}
			\left\|
			(\mathbf I-\mathbf P_{\bm h})
			\mathbf Q_{\bm h}
			\right\|_{\bm h,2}
			&\leq
			\big\|
			\big(\mathbf H_{\bm h}^{(\bm{k})}-\lambda^{(\bm{k})}\mathbf I\big)
			(\mathbf I-\mathbf P_{\bm h})
			\mathbf Q_{\bm h}
			\big\|_{\bm h,2}\\
			&\leq
			\big\|
			\big(\mathbf H_{\bm h}^{(\bm{k})}-\lambda^{(\bm{k})}\mathbf I\big)
			\mathbf Q_{\bm h}
			\big\|_{\bm h,2},
		\end{aligned}
		\]
		which proves \eqref{eq:err subspace gap bound}.
		
		As $\mathbf P_{\bm h}$ is a spectral projector of
		$\mathbf H_{\bm h}^{(\bm{k})}$, it commutes with the matrix:
		\[
		\mathbf H_{\bm h}^{(\bm{k})}\mathbf P_{\bm h}
		=
		\mathbf P_{\bm h}\mathbf H_{\bm h}^{(\bm{k})}.
		\]
		Hence,
		\begin{equation}\label{eq:err projected residual identity}
			\mathbf P_{\bm h}
			\big(
			\mathbf H_{\bm h}^{(\bm{k})}-\lambda^{(\bm{k})}\mathbf I
			\big)
			\mathbf P_{\bm h}
			\mathbf Q_{\bm h}=
			\mathbf P_{\bm h}
			\big(
			\mathbf H_{\bm h}^{(\bm{k})}-\lambda^{(\bm{k})}\mathbf I
			\big)
			\mathbf Q_{\bm h}.
		\end{equation}
		For sufficiently small $h_1$ and sufficiently large $N_2$,
		\eqref{eq:err subspace gap bound} gives
		$$\begin{aligned}
			\big\|\mathbf P_{\bm h}\mathbf Q_{\bm h}\widetilde{\bm W}\big\|_{\bm h}^2=&\big\|\mathbf Q_{\bm h}\widetilde{\bm W}\big\|_{\bm h}^2- \big\|(\mathbf I_{N_1N_2}-\mathbf P_{\bm h})\mathbf Q_{\bm h}\widetilde{\bm W}\big\|_{\bm h}^2\\
			\geq&\big(1-C(h_1^2+N_2^{-s})\big)\big\|\widetilde{\bm W}\big\|_{\bm h}^2,
		\end{aligned}$$
		Thus $\mathbf P_{\bm h}\mathbf Q_{\bm h}$ is an invertible map from $\mathbb C^m$
		onto $\mathcal M_{\bm h}(\lambda^{(\bm{k})})$, with
		$$\big\|(\mathbf P_{\bm h}\mathbf Q_{\bm h})^{-1}\big\|_{\bm h}\leq\big(1-C(h_1^2+N_2^{-s})\big)^{-\frac12}.$$
		Therefore, combining \eqref{eq:err Ah residual proof} with
		\eqref{eq:err projected residual identity} gives
		\begin{equation*}
			\big\|
			\mathbf P_{\bm h}
			\big(
			\mathbf H_{\bm h}^{(\bm{k})}-\lambda^{(\bm{k})}\mathbf I
			\big)
			\mathbf P_{\bm h}
			\big\|_{\bm h}
			\leq
			C\left(h_1^2+N_2^{-s}\right).
		\end{equation*}
		Note that $\mathbf P_{\bm h}$ is the discrete spectral projector associated with $I_{\delta_{\lambda^{(\bm{k})}}}$.
		Therefore, the Hermitian operator
		\[
		\mathbf P_{\bm h}
		\big(
		\mathbf H_{\bm h}^{(\bm{k})}-\lambda^{(\bm{k})}\mathbf I
		\big)
		\mathbf P_{\bm h}
		\]
		has eigenvalues
		\[
		\lambda_{i,\bm h}^{(\bm{k})}-\lambda^{(\bm{k})},
		\qquad i=1,\ldots,m,
		\]
		on $\mathcal M_{\bm h}(\lambda^{(\bm{k})})$, and vanishes on $\mathcal M_{\bm h}(\lambda^{(\bm{k})})^\perp$.
		It follows that
		\[
		\max_{1\leq i\leq m}
		\left|
		\lambda_{i,\bm h}^{(\bm{k})}-\lambda^{(\bm{k})}
		\right|
		=
		\big\|
		\mathbf P_{\bm h}
		\big(
		\mathbf H_{\bm h}^{(\bm{k})}-\lambda^{(\bm{k})}\mathbf I
		\big)
		\mathbf P_{\bm h}
		\big\|_{\bm h}\leq
		C\left(h_1^2+N_2^{-s}\right).
		\]
		which establishes \eqref{eq:err eigenvalue cluster bound}.
	\end{proof}
	
	For any two nonzero subspaces $\mathcal X_{\bm h}$ and $\mathcal Y_{\bm h}$
	of the discrete Fourier coefficient space, define the one-sided gap from
	$\mathcal X_{\bm h}$ to $\mathcal Y_{\bm h}$ by
	\begin{equation*}
		\begin{aligned}
			\operatorname{gap}_{\bm h}
			\left(
			\mathcal X_{\bm h},
			\mathcal Y_{\bm h}
			\right)
			&:=
			\sup_{\substack{
					\widetilde{\bm v}_{\bm h}\in\mathcal X_{\bm h}\\
					\|\widetilde{\bm v}_{\bm h}\|_{\bm h}=1}}
			\operatorname{dist}_{\bm h}
			\left(
			\widetilde{\bm v}_{\bm h},
			\mathcal Y_{\bm h}
			\right),
		\end{aligned}
	\end{equation*}
	where
	\[
	\operatorname{dist}_{\bm h}
	\left(
	\widetilde{\bm v}_{\bm h},
	\mathcal Y_{\bm h}
	\right)
	:=
	\inf_{\widetilde{\bm w}_{\bm h}\in\mathcal Y_{\bm h}}
	\left\|
	\widetilde{\bm v}_{\bm h}
	-
	\widetilde{\bm w}_{\bm h}
	\right\|_{\bm h}.
	\]
	The symmetric gap between the two subspaces is then defined by
	\begin{equation}\label{eq:err symmetric subspace gap definition}
		\begin{aligned}
			\widehat{\operatorname{gap}}_{\bm h}
			\left(
			\mathcal X_{\bm h},
			\mathcal Y_{\bm h}
			\right)
			&:=
			\max
			\left\{
			\operatorname{gap}_{\bm h}
			\left(
			\mathcal X_{\bm h},
			\mathcal Y_{\bm h}
			\right),
			\operatorname{gap}_{\bm h}
			\left(
			\mathcal Y_{\bm h},
			\mathcal X_{\bm h}
			\right)
			\right\}.
		\end{aligned}
	\end{equation}
	
	\begin{corollary}
		\label{cor:err full subspace gap}
		Under the assumptions of Theorem \ref{thm:err cluster}, let $\mathcal N_{\bm h}(\lambda^{(\bm{k})})$ be the subspace spanned by the columns of $\mathbf Q_{\bm h}$.
		Then $\mathcal N_{\bm h}(\lambda^{(\bm{k})})$ and the discrete eigenspace
		$\mathcal M_{\bm h}(\lambda^{(\bm{k})})$ satisfy
		\begin{equation}\label{eq:err full subspace gap bound}
			\widehat{\operatorname{gap}}_{\bm h}
			\left(
			\mathcal N_{\bm h}(\lambda^{(\bm{k})}),
			\mathcal M_{\bm h}(\lambda^{(\bm{k})})
			\right)
			\leq
			\frac{C}{\delta_{\lambda^{(\bm{k})}}}
			\left(
			h_1^2+N_2^{-s}
			\right).
		\end{equation}
	\end{corollary}
	
	\begin{proof}
		Since the columns of $\mathbf Q_{\bm h}$ are orthonormal with respect to the discrete inner product
		$\langle\cdot,\cdot\rangle_{\bm h}$,
		the map
		\[
		\mathbf Q_{\bm h}:
		\mathbb C^m\longrightarrow\mathcal N_{\bm h}(\lambda^{(\bm{k})})
		\]
		is an isometric isomorphism. Thus, as $\bm c$ ranges over
		all unit vectors in $\mathbb C^m$,
		$\mathbf Q_{\bm h}\bm c$ ranges over
		all unit vectors in $\mathcal N_{\bm h}(\lambda^{(\bm{k})})$.
		
		On the other hand, $\mathbf P_{\bm h}$ is the discrete orthogonal projector onto
		$\mathcal M_{\bm h}(\lambda^{(\bm{k})})$. Hence, for any
		$\widetilde{\bm v}_{\bm h}$,
		\[
		\operatorname{dist}_{\bm h}
		\left(
		\widetilde{\bm v}_{\bm h},
		\mathcal M_{\bm h}(\lambda^{(\bm{k})})
		\right)
		=
		\left\|
		\left(
		\mathbf I-\mathbf P_{\bm h}
		\right)
		\widetilde{\bm v}_{\bm h}
		\right\|_{\bm h}.
		\]
		Therefore,
		\begin{equation*}
			\begin{aligned}
				\operatorname{gap}_{\bm h}
				\left(
				\mathcal N_{\bm h}(\lambda^{(\bm{k})}),
				\mathcal M_{\bm h}(\lambda^{(\bm{k})})
				\right)
				&=
				\sup_{\substack{
						\widetilde{\bm v}_{\bm h}\in\mathcal N_{\bm h}(\lambda^{(\bm{k})})\\
						\|\widetilde{\bm v}_{\bm h}\|_{\bm h}=1}}
				\left\|
				\left(
				\mathbf I-\mathbf P_{\bm h}
				\right)
				\widetilde{\bm v}_{\bm h}
				\right\|_{\bm h}\\
				&=
				\sup_{\|\bm c\|_2=1}
				\left\|
				\left(
				\mathbf I-\mathbf P_{\bm h}
				\right)
				\mathbf Q_{\bm h}\bm c
				\right\|_{\bm h}\\
				&=
				\left\|
				\left(
				\mathbf I-\mathbf P_{\bm h}
				\right)
				\mathbf Q_{\bm h}
				\right\|_{\bm h,2}.
			\end{aligned}
		\end{equation*}
		By \eqref{eq:err subspace gap bound} in Theorem \ref{thm:err cluster},
		we have
		\begin{equation}\label{eq:err first directed gap bound}
			\operatorname{gap}_{\bm h}
			\left(
			\mathcal N_{\bm h}(\lambda^{(\bm{k})}),
			\mathcal M_{\bm h}(\lambda^{(\bm{k})})
			\right)
			\leq
			\frac{C}{\delta_{\lambda^{(\bm{k})}}}
			\left(
			h_1^2+N_2^{-s}
			\right).
		\end{equation}
		
		We next show that the one-sided gap in the reverse direction satisfies the same bound. Let
		$\mathbf P_{\mathcal N,\bm h}$ be the discrete orthogonal projector onto
		$\mathcal N_{\bm h}(\lambda^{(\bm{k})})$, and consider the restricted operator
		\[
		\mathbf S_{\bm h}
		:=
		\left.
		\mathbf P_{\bm h}
		\right|_{\mathcal N_{\bm h}(\lambda^{(\bm{k})})}
		:
		\mathcal N_{\bm h}(\lambda^{(\bm{k})})
		\longrightarrow
		\mathcal M_{\bm h}(\lambda^{(\bm{k})}).
		\]
		Its adjoint with respect to the discrete inner product is
		\[
		\mathbf S_{\bm h}^{*}
		=
		\left.
		\mathbf P_{\mathcal N,\bm h}
		\right|_{\mathcal M_{\bm h}(\lambda^{(\bm{k})})}
		:
		\mathcal M_{\bm h}(\lambda^{(\bm{k})})
		\longrightarrow
		\mathcal N_{\bm h}(\lambda^{(\bm{k})}).
		\]
		
		For any
		$\widetilde{\bm W}\in\mathcal N_{\bm h}(\lambda^{(\bm{k})})$,
		the orthogonal decomposition gives
		\[
		\big\|
		(\mathbf I-\mathbf P_{\bm h})
		\widetilde{\bm W}
		\big\|_{\bm h}^2
		=
		\big\|\widetilde{\bm W}\big\|_{\bm h}^2
		-
		\big\|\mathbf P_{\bm h}
		\widetilde{\bm W}\big\|_{\bm h}^2
		=
		\big\|\widetilde{\bm W}\big\|_{\bm h}^2
		-
		\big\|\mathbf S_{\bm h}
		\widetilde{\bm W}\big\|_{\bm h}^2.
		\]
		Hence,
		\begin{equation}\label{eq:err directed gap singular value one}
			\begin{aligned}
				\operatorname{gap}_{\bm h}^2
				\left(
				\mathcal N_{\bm h}(\lambda^{(\bm{k})}),
				\mathcal M_{\bm h}(\lambda^{(\bm{k})})
				\right)
				=
				1-\sigma_{\min}^2(\mathbf S_{\bm h}),
			\end{aligned}
		\end{equation}
		where
		$$\sigma_{\min}(\mathbf S_{\bm h})=\inf_{
			\|\widetilde{\bm W}\|_{\bm h}=1}
		\big\|
		\mathbf S_{\bm h}
		\widetilde{\bm W}
		\big\|_{\bm h}$$
		denotes
		the smallest singular value of $\mathbf S_{\bm h}$. Moreover, $\sigma_{\min}^2(\mathbf S_{\bm h})$ is the smallest eigenvalue of the Hermitian operator $\mathbf S_{\bm h}^*\mathbf S_{\bm h}$.
		Similarly,
		\begin{equation}\label{eq:err directed gap singular value two}
			\begin{aligned}
				\operatorname{gap}_{\bm h}^2
				\left(
				\mathcal M_{\bm h}(\lambda^{(\bm{k})}),
				\mathcal N_{\bm h}(\lambda^{(\bm{k})})
				\right)
				=
				1-
				\sigma_{\min}^2
				\left(
				\mathbf S_{\bm h}^{*}
				\right).
			\end{aligned}
		\end{equation}
		By the dimension-preservation assumption for the discrete eigenvalue cluster in Theorem \ref{thm:err cluster},
		\[
		\dim\mathcal N_{\bm h}(\lambda^{(\bm{k})})
		=
		\dim\mathcal M_{\bm h}(\lambda^{(\bm{k})})
		=
		m.
		\]
		Thus $\mathbf S_{\bm h}$ maps between two finite-dimensional inner-product spaces of the same dimension,
		and $\mathbf S_{\bm h}$ and
		$\mathbf S_{\bm h}^{*}$ have the same smallest singular value. Therefore,
		\[
		\sigma_{\min}(\mathbf S_{\bm h})
		=
		\sigma_{\min}
		\left(
		\mathbf S_{\bm h}^{*}
		\right).
		\]
		Combining \eqref{eq:err directed gap singular value one} and
		\eqref{eq:err directed gap singular value two} yields
		\begin{equation*}
			\operatorname{gap}_{\bm h}
			\left(
			\mathcal N_{\bm h}(\lambda^{(\bm{k})}),
			\mathcal M_{\bm h}(\lambda^{(\bm{k})})
			\right)
			=
			\operatorname{gap}_{\bm h}
			\left(
			\mathcal M_{\bm h}(\lambda^{(\bm{k})}),
			\mathcal N_{\bm h}(\lambda^{(\bm{k})})
			\right).
		\end{equation*}
		By the definition \eqref{eq:err symmetric subspace gap definition} of the symmetric gap
		and the bound \eqref{eq:err first directed gap bound}, we obtain
		\eqref{eq:err full subspace gap bound}.
	\end{proof}
	
	\begin{remark}
		\label{rem:cluster-dimension-source-resolvent}
		Assumption \ref{ass:err discrete cluster dimension} excludes additional spurious eigenvalues from the target interval, and the eigenvalue and eigenspace estimates in Theorem \ref{thm:err cluster} are conditional on this assumption. Establishing dimension preservation requires additional spectral convergence analysis; see \cite{BabuskaOsborn1991,DNR1978a,CancesEhrlacherMaday2014,Gallistl2015} for related results on self-adjoint eigenproblems and eigenvalue clusters.
	\end{remark}
	
	For a simple isolated eigenvalue, the preceding eigenspace error estimate further yields an error bound for the discrete eigenvector after an appropriate phase adjustment.
	
	\begin{corollary}
		\label{cor:err simple eigen}
		Suppose $\lambda^{(\bm{k})}$ is a simple isolated eigenvalue, and let $u$ be an associated eigenfunction with
		$\|u\|_{L^2(\Omega)}=1$. Set
		\[
		\widetilde{\bm z}_{\bm h}
		:=
		\frac{
			\mathbf{\widetilde R}_{\bm h}u
		}{
			\|\mathbf{\widetilde R}_{\bm h}u\|_{\bm h}
		}.
		\]
		Let $\lambda_{\bm h}^{(\bm{k})}$ be the discrete eigenvalue near $\lambda^{(\bm{k})}$,
		and let $\widetilde{\bm U}_{\bm h}$ be an associated discrete eigenvector with
		$\|\widetilde{\bm U}_{\bm h}\|_{\bm h}=1$.
		Then there exists a phase $\theta_{\bm h}\in\mathbb R$ such that
		\begin{equation}\label{eq:err simple eigenvector}
			\left\|
			\widetilde{\bm U}_{\bm h}
			-
			\mathrm e^{\mathrm i\theta_{\bm h}}
			\widetilde{\bm z}_{\bm h}
			\right\|_{\bm h}
			\leq
			\frac{C}{\delta_{\lambda^{(\bm{k})}}}
			\left(h_1^2+N_2^{-s}\right),
		\end{equation}
		Furthermore,
		\begin{equation*}
			|\lambda_{\bm h}^{(\bm{k})}-\lambda^{(\bm{k})}|
			\leq
			C\left(h_1^2+N_2^{-s}\right).
		\end{equation*}
		By the Parseval identity \eqref{eq:err parseval physical Fourier},
		the estimate \eqref{eq:err simple eigenvector} is equivalent to a discrete
		$L^2$-type error estimate at the physical grid points.
	\end{corollary}
	
	\vskip 0.8cm

	\section{Numerical Experiments}\setcounter{equation}{0}

	Three numerical examples assess the proposed hybrid Fourier--finite difference scheme in terms of band functions, eigenstate densities, and errors in eigenvalues, eigenstates, and eigenspaces. The discrete eigenvalue problems are solved by the \gls{LOBPCG} method, and the preconditioner is described in Section 3.3.
	
	For each quasimomentum \(\bm k\), let \(\lambda_j^{(\bm k)}\), \(\lambda_{j,\bm h}^{(\bm k)}\), and \(\lambda_{j,\mathrm{ref}}^{(\bm k)}\) denote the \(j\)th eigenvalues of the continuous problem, the computed discretization, and the high-resolution reference discretization, respectively, ordered nondecreasingly and counted with multiplicity. The computed and reference eigenpairs use the same model and magnetic unit cell but different grids. Their eigenvectors are stored as Fourier coefficients in the \(y\)-direction and normalized in the discrete inner product of Section 3.1. Band diagrams plot \(\lambda_{j,\bm h}^{(\bm k)}\), while \(u_j\) in density plots denotes the corresponding eigenstate in physical space.
	
	\subsection{Magnetic field and potential with the same spatial modulation}
	Let $a=b=2$, and choose the magnetic field and scalar potential as
	\begin{equation}\label{eq:dirac model}
		B(x,y)=\pi+\exp\big(\cos\theta(x,y)\big)-I_0,
	\end{equation}
	\begin{equation}\label{eq:dirac V}
		V(x,y)=\exp\big(\cos\theta(x,y)\big),
	\end{equation}
	where $\theta(x,y)=\pi x+\pi y$ and
	$$I_0=\frac{1}{2\pi}\int_0^{2\pi}\mathrm e^{\cos t}\,\mathrm dt.$$
	It follows that $\bar B=\pi$ and $n_B=2$, so the integer magnetic flux condition \eqref{eq:integer flux} is satisfied.
	
	We compute the first six band functions over the rectangular Brillouin zone
	\begin{equation*}
		\mathcal B=\left[0,\frac{2\pi}{a}\right]\times
		\left[0,\frac{2\pi}{b}\right]
	\end{equation*}
	The four vertices are denoted, in order, by $P_1=(0,0)$, $P_2=(2\pi/a,0)$, $P_3=(2\pi/a,2\pi/b)$, and $P_4=(0,2\pi/b)$. Dividing the two quasimomentum intervals into $m$ and $n$ equal subintervals, respectively, gives the sampling set including the boundary endpoints:
	\begin{equation*}
		\mathcal B_{m\times n}
		=\left\{\left(\frac{2\pi\ell_1}{am},\frac{2\pi\ell_2}{bn}\right):
		\ell_1=0,\ldots,m,\ \ell_2=0,\ldots,n\right\}.
	\end{equation*}
	Here $m,n\in\mathbb N_+$ are the numbers of subdivisions in the respective directions. The stopping tolerance for \gls{LOBPCG} in this example is set to $10^{-8}$.
	
\begin{figure}[H]\centering
		\includegraphics[width=0.72\textwidth]{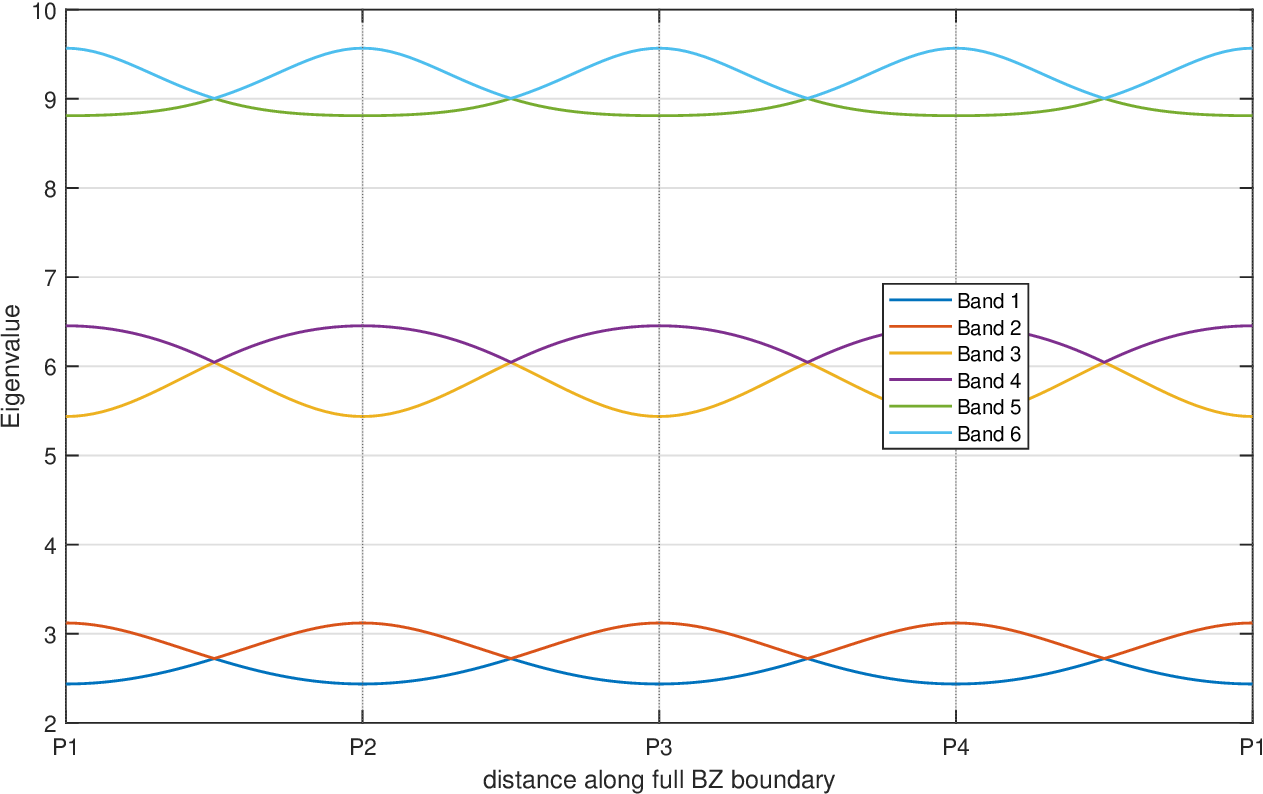}
		\caption{The first six numerical band functions along the path $P_1$--$P_2$--$P_3$--$P_4$--$P_1$.}\label{fig:dirac bands}
	\end{figure}

\begin{figure}[H]\centering
		\includegraphics[width=0.96\textwidth]{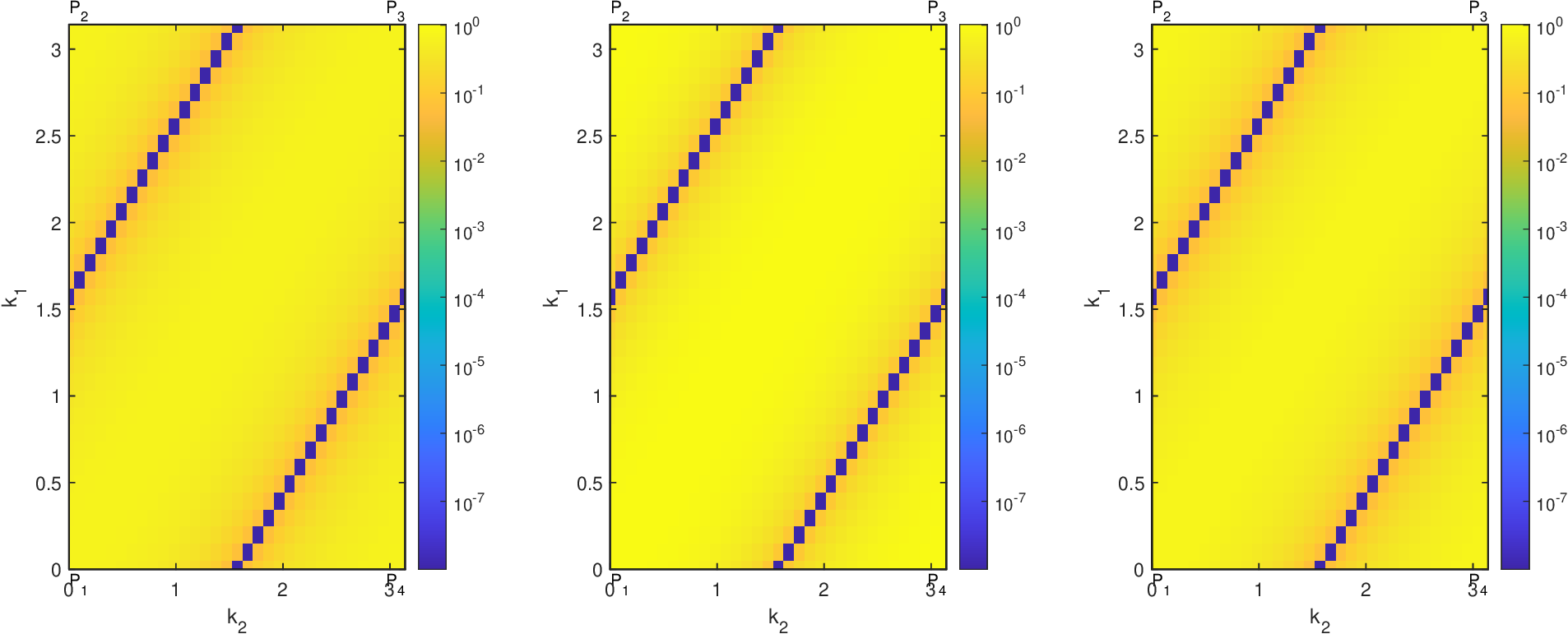}
		\caption{Pairwise eigenvalue gaps $g_{12}$, $g_{34}$, and $g_{56}$ on the sampling set $\mathcal B_{32\times32}$.}\label{fig:dirac pair gaps}
	\end{figure}

	Figure~\ref{fig:dirac bands} shows the first six numerical band functions along the closed path $P_1$--$P_2$--$P_3$--$P_4$--$P_1$. Using the numerical reference eigenvalues, we define the pairwise eigenvalue gaps by
	\begin{equation*}
		g_{2j-1,2j}(\bm k)
		=\left|\lambda_{2j,\mathrm{ref}}^{(\bm k)}
		-\lambda_{2j-1,\mathrm{ref}}^{(\bm k)}\right|,
		\qquad j=1,2,3.
	\end{equation*}
	Figure~\ref{fig:dirac pair gaps} displays the distributions of $g_{12}$, $g_{34}$, and $g_{56}$ over $\mathcal B_{32\times32}$. The nearly vanishing pairwise gaps lie along diagonal lines passing through the midpoints of the boundary edges; under periodic identification of the boundaries, these lines form closed structures.
	
\begin{figure}[H]\centering
		\includegraphics[width=1\textwidth]{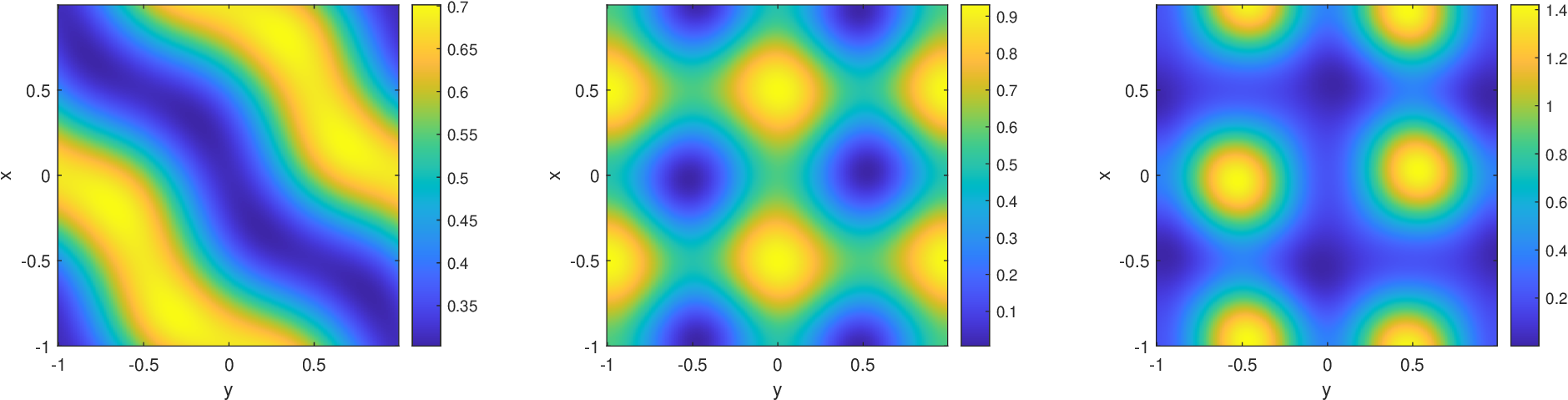}
		\caption{Total densities of the three two-dimensional degenerate eigenspaces at $\bm k=(\pi/a,0)$.}\label{fig:dirac subspace density}
	\end{figure}

	At $\bm k=(0,0)$, the first six eigenvalues are distinct, whereas at $\bm k=(\pi/a,0)$ the eigenvalue pairs $1$--$2$, $3$--$4$, and $5$--$6$ are degenerate. To avoid the dependence of individual-state densities on the choice of basis within a degenerate eigenspace, we consider the total densities
	\begin{equation*}
		\rho_{12}=|u_1|^2+|u_2|^2,\qquad
		\rho_{34}=|u_3|^2+|u_4|^2,\qquad
		\rho_{56}=|u_5|^2+|u_6|^2,
	\end{equation*}
	to characterize the spatial distributions of the corresponding two-dimensional eigenspaces. The results are shown in Figure~\ref{fig:dirac subspace density}.
	
	Taking the numerical solution computed on the $2048\times256$ grid as the reference solution, we define the relative error in the $j$th eigenvalue by
	\begin{equation*}
		e_{\mathrm{eig}}^{(j)}(\bm k)
		=\frac{\left|\lambda_{j,\bm h}^{(\bm k)}-\lambda_{j,\mathrm{ref}}^{(\bm k)}\right|}
		{\left|\lambda_{j,\mathrm{ref}}^{(\bm k)}\right|}.
	\end{equation*}
	All errors reported below are measured relative to the specified numerical reference solutions and should not be interpreted as errors relative to the exact solutions of the continuous problem. Whenever a maximum over quasimomenta is reported, it is taken only over the discrete sampling points used in the corresponding experiment.
	
	Eigenstate errors are computed directly from the stored Fourier coefficients. In the following inner products and matrix products, the coefficient vectors have consistent orderings of the $x$-grid nodes and Fourier-mode indices, and the discrete inner product of Section~3.1 is used. Let $\bm Q_{\mathrm{ref}}$ and $\bm Q_{\mathrm{test}}$ be matrices whose columns form orthonormal bases of the reference and computed subspaces, respectively, represented by Fourier coefficient vectors. Thus, $h_1b\bm Q_{\mathrm{ref}}^*\bm Q_{\mathrm{ref}}=h_1b\bm Q_{\mathrm{test}}^*\bm Q_{\mathrm{test}}=\mathbf I$. Their cross-Gram matrix is
	\begin{equation}\label{eq:esub gram}
		\mathbf G_{\mathrm{cross},\bm h}
		=h_1b\bm Q_{\mathrm{ref}}^*\bm Q_{\mathrm{test}},
	\end{equation}
	The matrix $\mathbf G_{\mathrm{cross},\bm h}$ is distinct from the sampling Gram matrix $\mathbf G_{\bm h}$ introduced in Section~4. For brevity, the dependence of the basis matrices and the cross-Gram matrix on $\bm k$ is suppressed. We define the subspace error by
	\begin{equation}\label{eq:esub}
		e_{\mathrm{sub}}(\bm k)
		=\sqrt{\max\left\{0,1-
			\sigma_{\min}\big(\mathbf G_{\mathrm{cross},\bm h}\big)^2\right\}}.
	\end{equation}
	This quantity equals the sine of the largest principal angle between the two subspaces and is independent of the choice of their orthonormal bases. For the subspace spanned by eigenstates $j_1$ through $j_2$, we write the corresponding error as $e_{\mathrm{sub}}^{j_1:j_2}(\bm k)$; the superscript $j_1:j_2$ denotes the inclusive range of consecutive eigenstate indices.
	
	For the $j$th eigenstate associated with a simple eigenvalue, let $\widetilde{\bm U}_{j,\bm h}^{(\bm k)}$ and $\widetilde{\bm U}_{j,\mathrm{ref},\bm h}^{(\bm k)}$ denote the computed and reference Fourier coefficient vectors, respectively, both normalized in the discrete inner product above. We define the optimally phase-aligned error by
	\begin{equation}\label{eq:single eu}
		\begin{aligned}
			e_u^{(j)}(\bm k)
			=\min_{\theta_{\bm h}\in\mathbb R}
			\left\|\widetilde{\bm U}_{j,\bm h}^{(\bm k)}
			-\mathrm e^{\mathrm i\theta_{\bm h}}
			\widetilde{\bm U}_{j,\mathrm{ref},\bm h}^{(\bm k)}\right\|_{\bm h}
			=\sqrt{2-2\left|
				\left\langle\widetilde{\bm U}_{j,\mathrm{ref},\bm h}^{(\bm k)},
				\widetilde{\bm U}_{j,\bm h}^{(\bm k)}\right\rangle_{\bm h}\right|}.
		\end{aligned}
	\end{equation}
	The parameter $\theta_{\bm h}$ has the same interpretation as the global phase in Corollary~\ref{cor:err simple eigen}.
	
	Fix $N_2=256$ and take $N_1=128,256,512$. Table~\ref{tab:dirac x convergence} reports the maximum relative errors in the first six eigenvalues over the quasimomentum sampling points, together with the error in the first eigenstate at $\bm k=(0,0)$. Table~\ref{tab:dirac_d1_subspaces} reports the errors in the three two-dimensional degenerate eigenspaces and the combined subspace spanned by the first six eigenstates at $\bm k=(\pi/a,0)$. For successive doublings of the resolution in the $x$-direction, the observed convergence order is defined by
	\begin{equation*}
		p_{\mathrm{obs}}(N_1\to2N_1)
		=\log_2\frac{e(N_1)}{e(2N_1)},
	\end{equation*}
	where $e(N_1)$ denotes the corresponding error measured against the same numerical reference solution. The observed orders in both tables are close to $2$, consistent with the second-order error estimates in Section~4. Table~\ref{tab:dirac y convergence} presents the numerical errors for different numbers of Fourier modes with $N_1=1024$ fixed. The errors decrease rapidly as $N_2$ increases and then level off, reflecting the high-order accuracy of the Fourier discretization in the $y$-direction until the remaining error is dominated by the $x$-direction discretization.
	
	\begin{table}[htbp]\centering\small\setlength{\tabcolsep}{6pt}
		\caption{Eigenvalue and individual-eigenstate errors with $N_2=256$ fixed.}\label{tab:dirac x convergence}
		\begin{tabular}{ccc}
			\hline
			$N_1$ & $\max\limits_{\bm k\in \mathcal B_{32\times32},\,1\le j\le6}e_{\mathrm{eig}}^{(j)}(\bm k)$ & $e_u^{(1)}(0,0)$\\
			\hline
			$128$ & $1.1\times10^{-4}$ & $2.2\times10^{-5}$\\
			$256$ & $2.7\times10^{-5}$ & $5.4\times10^{-6}$\\
			$512$ & $6.5\times10^{-6}$ & $1.3\times10^{-6}$\\
			\hline
			$p_{\mathrm{obs}}(128\to256)$ & $2.01723$ & $2.01750$\\
			$p_{\mathrm{obs}}(256\to512)$ & $2.07043$ & $2.07366$\\
			\hline
		\end{tabular}
	\end{table}
	
	\begin{table}[htbp]\centering\small\setlength{\tabcolsep}{5pt}
		\caption{Eigenspace errors with $N_2=256$ fixed.}\label{tab:dirac_d1_subspaces}
		\begin{tabular}{ccccc}
			\hline
			$N_1$ & $e_{\mathrm{sub}}^{1:2}(\pi/a,0)$ & $e_{\mathrm{sub}}^{3:4}(\pi/a,0)$ & $e_{\mathrm{sub}}^{5:6}(\pi/a,0)$ & $e_{\mathrm{sub}}^{1:6}(\pi/a,0)$\\
			\hline
			$128$ & $3.70\times10^{-5}$ & $1.14\times10^{-4}$ & $2.07\times10^{-4}$ & $2.05\times10^{-4}$\\
			$256$ & $9.14\times10^{-6}$ & $2.81\times10^{-5}$ & $5.12\times10^{-5}$ & $5.07\times10^{-5}$\\
			$512$ & $2.18\times10^{-6}$ & $6.69\times10^{-6}$ & $1.22\times10^{-5}$ & $1.21\times10^{-5}$\\
			\hline
			$p_{\mathrm{obs}}(128\to256)$ & $2.01720$ & $2.01734$ & $2.01740$ & $2.01741$\\
			$p_{\mathrm{obs}}(256\to512)$ & $2.07048$ & $2.07046$ & $2.07047$ & $2.07047$\\
			\hline
		\end{tabular}
	\end{table}

	\begin{table}[htbp]\centering\small\setlength{\tabcolsep}{8pt}
		\caption{Eigenvalue and eigenspace errors with $N_1=1024$ fixed.}\label{tab:dirac y convergence}
		\begin{tabular}{ccc}
			\hline
			$N_2$ & $\max\limits_{\bm k\in \mathcal B_{32\times32},\,1\le j\le6}e_{\mathrm{eig}}^{(j)}(\bm k)$ & $\max\limits_{\bm k\in \mathcal B_{32\times32}}e_{\mathrm{sub}}^{1:6}(\bm k)$\\
			\hline
			$4$  & $1.2301\times10^{-1}$ & $5.2514\times10^{-1}$\\
			$8$  & $2.3079\times10^{-6}$ & $8.9222\times10^{-4}$\\
			$12$ & $1.3057\times10^{-6}$ & $4.0667\times10^{-6}$\\
			$16$ & $1.3057\times10^{-6}$ & $2.7040\times10^{-6}$\\
			$20$ & $1.3057\times10^{-6}$ & $2.7040\times10^{-6}$\\
			$24$ & $1.3057\times10^{-6}$ & $2.7040\times10^{-6}$\\
			\hline
		\end{tabular}
	\end{table}
	
	\subsection{Periodic magnetic field with a multifrequency potential}
	
	This example considers a nonuniform periodic magnetic field and a multifrequency potential containing incommensurate spatial frequencies. Rational approximations of the frequencies are used to construct periodic test problems, for which we examine the spatial distributions of eigenstates and the accuracy of eigenvalue computations on enlarged magnetic unit cells. The magnetic field is given by
	\begin{equation}\label{eq:qp B}
		B(x,y)=\frac{2}{\pi}\big(1+\cos2x+\cos2y\big),
	\end{equation}
	and the scalar potential by
	\begin{equation}\label{eq:qp V}
		V(x,y)=2\big(\sin^2x+\sin^2y-
		\cos(2x+2\sqrt2\,y)\big).
	\end{equation}
	The potential contains an incommensurate spatial frequency in the $y$-direction. To obtain periodic test problems compatible with the magnetic Bloch discretization developed here, we approximate $\sqrt2$ by $p/q$ using rational approximation of the incommensurate spatial frequencies \cite{Jiang-Li-Zhang2025}, where $p,q$ are coprime positive integers. The quantity being approximated is the irrational spatial frequency in the potential, not the magnetic flux. Specifically, we take
	\begin{equation}\label{eq:qp approximations}
		\frac pq\in\left\{\frac32,\frac75,\frac{17}{12}\right\}.
	\end{equation}
	The corresponding periodic approximations of the potential are
	\begin{equation}\label{eq:qp V approx}
		V_{p,q}(x,y)=2\left(\sin^2x+\sin^2y-
		\cos\left(2x+\frac {2p}q\,y\right)\right).
	\end{equation}
	For each fixed $p/q$, we choose the magnetic unit cell with $a=\pi$ and $b=q\pi$. Both $B$ and $V_{p,q}$ are then periodic on this cell, and $\bar B=2/\pi$ and $n_B=q$, satisfying the integer magnetic flux condition \eqref{eq:integer flux}.

	To quantify the spatial concentration of an eigenstate, we recover the physical-grid vector $\bm U$ from a normalized Fourier coefficient vector $\|\widetilde{\bm U}\|_{\bm h}=1$ by the inverse \gls{DFT}, and define the \gls{IPR} by
	\begin{equation}\label{eq:qp IPR}
		\mathrm{IPR}
		=h_1h_2\sum_{l=0}^{N_1-1}\sum_{m=0}^{N_2-1}|U_{l,m}|^4.
	\end{equation}
	A larger \gls{IPR} generally indicates a more spatially concentrated probability density.
	
\begin{figure}[htbp]\centering
		\includegraphics[width=0.86\textwidth]{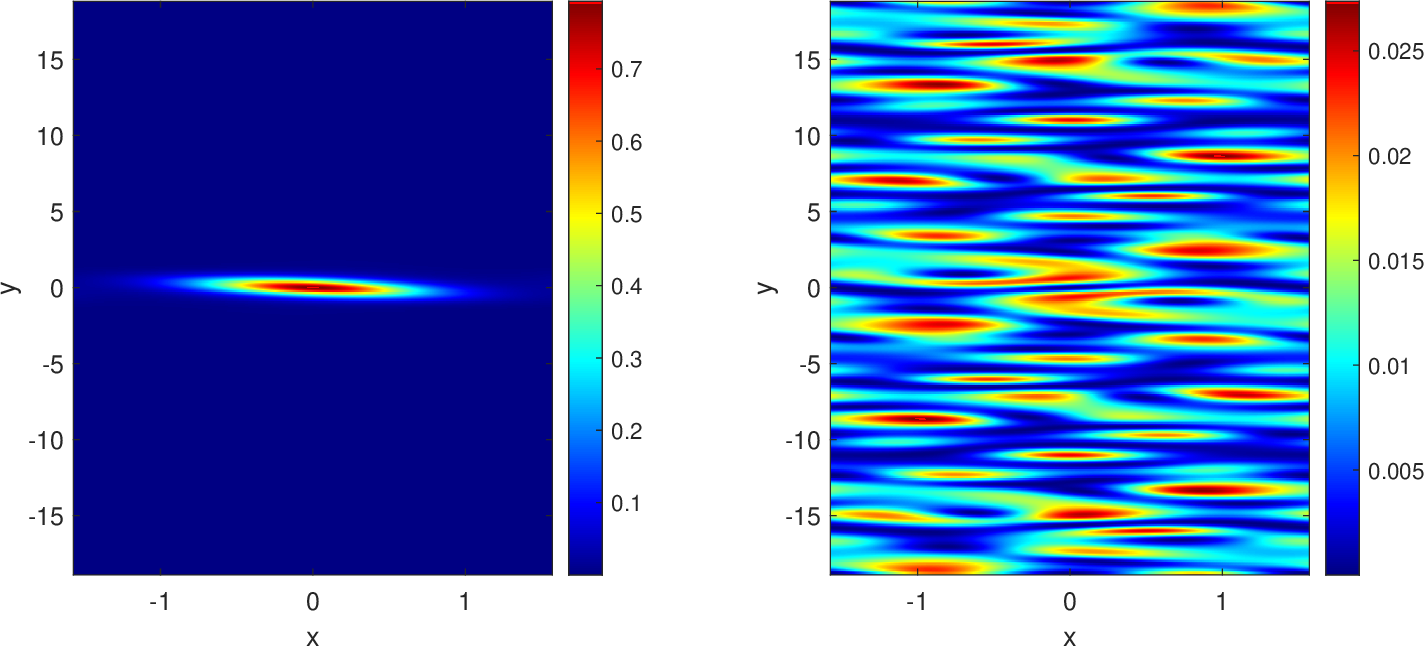}
		\caption{Probability densities of the eigenstates attaining the largest and smallest \gls{IPR} among the first $100$ eigenstates at $p/q=17/12$ and $\bm k=(0,0)$ (states $1$ and $61$, respectively).}\label{fig:qp fixed color ipr extrema}
	\end{figure}

	For $p/q=17/12$ and $\bm k=(0,0)$, we compute the first $100$ eigenstates. Among them, the first eigenstate has the largest \gls{IPR} $(3.65901\times10^{-1})$, while the $61$st has the smallest \gls{IPR} $(1.25698\times10^{-2})$. Their probability densities are displayed in Figure~\ref{fig:qp fixed color ipr extrema}. The probability density of the first eigenstate is concentrated mainly near $y=0$, whereas that of the $61$st eigenstate is distributed over a considerably wider region.
	
	We next examine grid errors in the low-lying eigenvalues at $\bm k=(0,0)$. For each of the three periodic problems obtained with $p/q=3/2,7/5,17/12$, we fix the number of Fourier modes at $N_2=256$, take $N_1=128,256,512,1024$, and use the solution on the $4096\times512$ grid as the numerical reference. Table~\ref{tab:qp eigenvalue errors} reports the maximum relative error among the first $100$ eigenvalues and the observed orders between successive grids. As $N_1$ doubles, the errors decrease consistently for all three approximations, with observed orders between $2.00$ and $2.07$, in agreement with the expected second-order behavior of the covariant central differences in the $x$-direction.
	
	\begin{table}[htbp]\centering\small\setlength{\tabcolsep}{8pt}
		\caption{Maximum relative error in the first $100$ eigenvalues for three enlarged magnetic unit cells with $N_2=256$ fixed; all reference solutions use a $4096\times512$ grid.}\label{tab:qp eigenvalue errors}
		\begin{tabular}{cccc}
			\hline
			$p/q$ & $N_1$ & $\max\limits_{1\le j\le100}e_{\mathrm{eig}}^{(j)}(0,0)$ & $p_{\mathrm{obs}}$\\
			\hline
			$3/2$ & $128$ & $1.45645\times10^{-3}$ & --\\
			$3/2$ & $256$ & $3.62959\times10^{-4}$ & $2.00457$\\
			$3/2$ & $512$ & $8.96645\times10^{-5}$ & $2.01720$\\
			$3/2$ & $1024$ & $2.13482\times10^{-5}$ & $2.07042$\\
			\hline
			$7/5$ & $128$ & $4.46854\times10^{-4}$ & --\\
			$7/5$ & $256$ & $1.11371\times10^{-4}$ & $2.00444$\\
			$7/5$ & $512$ & $2.75141\times10^{-5}$ & $2.01713$\\
			$7/5$ & $1024$ & $6.55092\times10^{-6}$ & $2.07040$\\
			\hline
			$17/12$ & $128$ & $1.95224\times10^{-4}$ & --\\
			$17/12$ & $256$ & $4.86618\times10^{-5}$ & $2.00427$\\
			$17/12$ & $512$ & $1.20223\times10^{-5}$ & $2.01708$\\
			$17/12$ & $1024$ & $2.86244\times10^{-6}$ & $2.07039$\\
			\hline
		\end{tabular}
	\end{table}

	\subsection{Periodic potential with a multifrequency magnetic field}\label{sec:qp54}
	
	We now consider a magnetic field containing components with different spatial frequencies in both coordinate directions. Let
	\begin{equation}\label{eq:qp54 original field}
		B(x,y)=\frac{2}{\pi}+\cos2x+\cos2y
		+\cos(2\sqrt2\,x)+\cos(2\sqrt2\,y),
	\end{equation}
	\begin{equation}\label{eq:qp54 potential}
		V(x,y)=2(\cos2x+\cos2y).
	\end{equation}
	The original magnetic field contains two incommensurate spatial frequencies in each of the $x$- and $y$-directions and therefore has no finite common period. As in the preceding example, we approximate the irrational spatial frequency in the magnetic field by replacing $\sqrt2$ with $p/q$, thereby obtaining a periodic magnetic field $B_{p,q}$ for numerical computation:
	\begin{equation}\label{eq:qp54 periodic field}
		B_{p,q}(x,y)=\frac{2}{\pi}+\cos2x+\cos2y
		+\cos\!\left(\frac{2p}{q}x\right)+\cos\!\left(\frac{2p}{q}y\right).
	\end{equation}
	Choosing the common magnetic unit cell with $a=b=q\pi$ gives
	\begin{equation}\label{eq:qp54 flux}
		\bar B=\frac{2}{\pi},\qquad
		n_B=\frac{\bar B ab}{2\pi}=q^2.
	\end{equation}
	In the computations below, we set $p/q=17/12$, corresponding to $a=b=12\pi$ and normalized magnetic flux $n_B=144$, and fix the quasimomentum at $\bm k=(0,0)$.
	
\begin{figure}[htbp]\centering
		
		\includegraphics[width=0.82\textwidth]{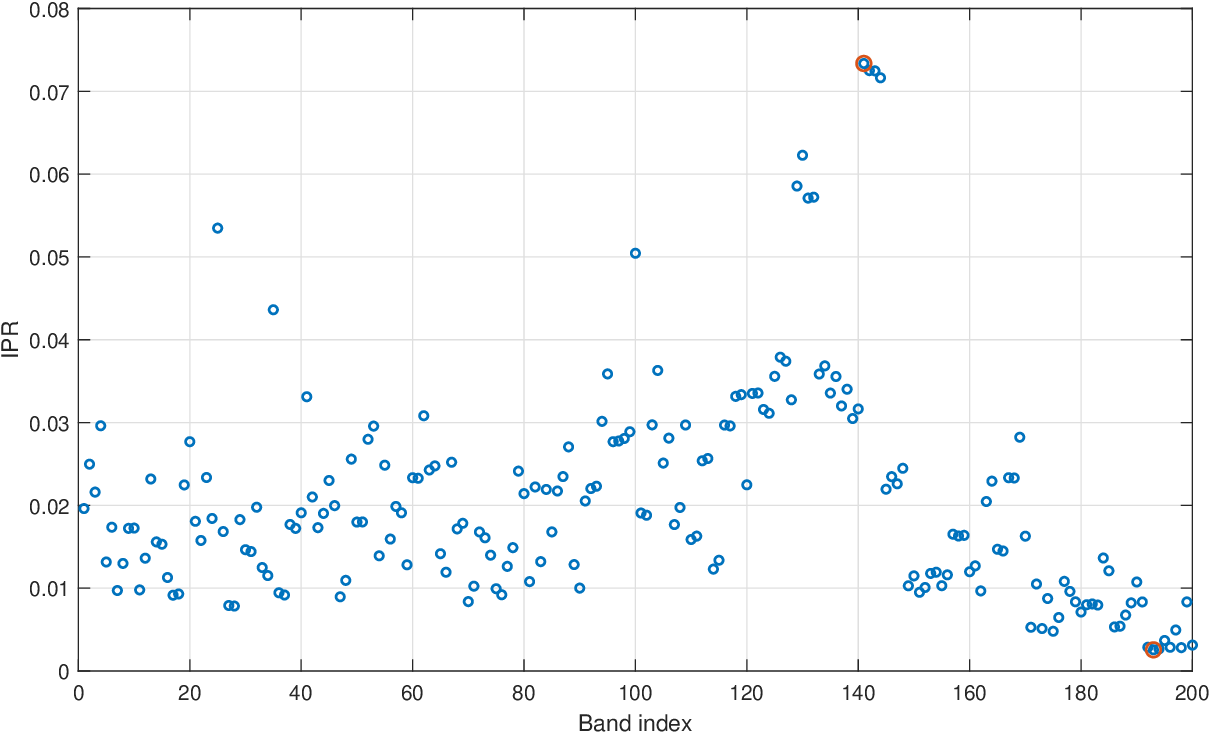}
		\caption{\gls{IPR} values of the first $200$ eigenstates.}\label{fig:qp54 ipr250}
	\end{figure}

\begin{figure}[htbp]\centering
		\includegraphics[width=0.98\textwidth]{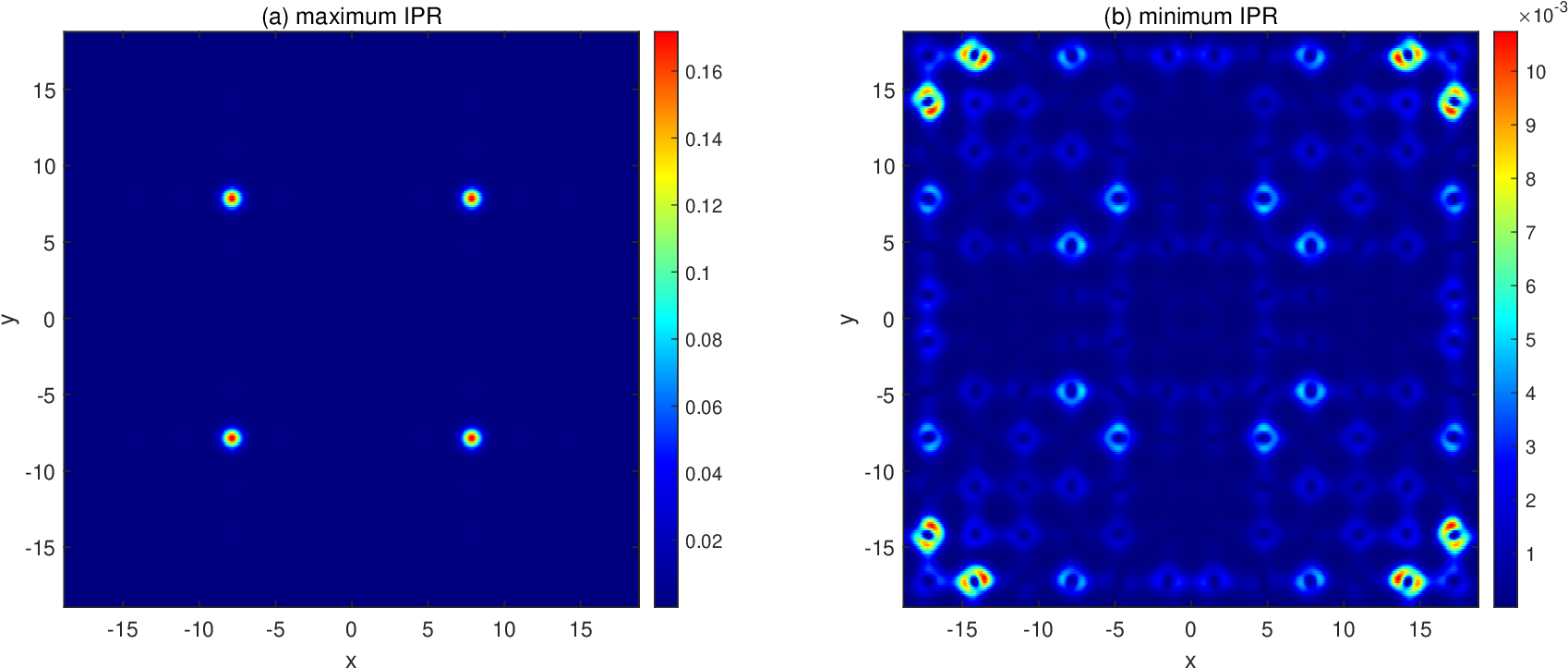}
		\caption{Probability densities of the eigenstates with extremal \gls{IPR} among the first $200$ states: left, the $141$st eigenstate with the largest \gls{IPR} in the computed eigenbasis; right, the $193$rd eigenstate with the smallest \gls{IPR}.}\label{fig:qp53 fixed color ipr extrema}
	\end{figure}

	We compute the first $200$ eigenpairs and use \eqref{eq:qp IPR} to measure the spatial concentration of the corresponding eigenstates. Figure~\ref{fig:qp54 ipr250} displays the \gls{IPR} values ordered by increasing eigenvalue; the orange open circles mark the maximum and minimum \gls{IPR} values. For the orthonormal set of computed eigenvectors in this run, the largest \gls{IPR}, $7.33639\times10^{-2}$, occurs for the $141$st eigenstate, whereas the smallest, $2.56012\times10^{-3}$, occurs for the $193$rd. Figure~\ref{fig:qp53 fixed color ipr extrema} shows their probability densities. The density of the $141$st eigenstate is concentrated in a few localized regions, while that of the $193$rd extends over a wider area and exhibits several local peaks.
	
	We further assess grid errors in the low-lying eigenvalues for fixed $p/q=17/12$ and $\bm k=(0,0)$. First, we fix $N_2=256$, take $N_1=128,256,512,1024$, and use the solution on the $4096\times1024$ grid as the numerical reference. Table~\ref{tab:qp54 x resolution} reports the maximum relative error among the first $20$ eigenvalues and the observed orders between successive grids. As $N_1$ increases, the maximum relative error decreases from $1.14739\times10^{-2}$ to $1.64664\times10^{-4}$. The observed orders for the three refinements are $2.02805$, $2.02283$, and $2.07182$, respectively, consistent with the expected second-order accuracy of the covariant discretization in the $x$-direction.
	
	To assess the effect of Fourier resolution, we next fix $N_1=1024$ and take $N_2=64,128,256,512$, retaining the same $4096\times1024$ reference grid. Table~\ref{tab:qp54 y resolution} shows a substantial reduction in the maximum relative error among the first $20$ eigenvalues as $N_2$ increases from $64$ to $256$. Increasing $N_2$ further to $512$ does not change the reported error at the displayed precision. A direct comparison between the $1024\times256$ and $1024\times512$ grids gives a maximum relative eigenvalue difference of only $4.07589\times10^{-11}$. Thus, at the present resolution, adding Fourier modes has a negligible effect on these $20$ eigenvalues; the remaining error relative to the high-resolution reference is primarily attributable to the $x$-direction discretization.
	
	\begin{table}[htbp]\centering\small\setlength{\tabcolsep}{9pt}
		\caption{Maximum relative error in the first $20$ eigenvalues and observed convergence orders with $N_2=256$ fixed.}\label{tab:qp54 x resolution}
		\begin{tabular}{ccc}
			\hline
			$N_1$ & $\max\limits_{1\le j\le20}e_{\mathrm{eig}}^{(j)}(0,0)$ & $p_{\mathrm{obs}}$\\
			\hline
			$128$ & $1.14739\times10^{-2}$ & --\\
			$256$ & $2.81326\times10^{-3}$ & $2.02805$\\
			$512$ & $6.92274\times10^{-4}$ & $2.02283$\\
			$1024$ & $1.64664\times10^{-4}$ & $2.07182$\\
			\hline
		\end{tabular}
	\end{table}
	
	\begin{table}[htbp]\centering\small\setlength{\tabcolsep}{12pt}
		\caption{Relative errors in the first $20$ eigenvalues with $N_1=1024$ fixed.}\label{tab:qp54 y resolution}
		\begin{tabular}{cc}
			\hline
			$N_2$ & $\max\limits_{1\le j\le20}e_{\mathrm{eig}}^{(j)}(0,0)$\\
			\hline
			$64$  & $4.32324\times10^{-2}$\\
			$128$ & $9.50210\times10^{-3}$\\
			$256$ & $1.64664\times10^{-4}$\\
			$512$ & $1.64664\times10^{-4}$\\
			\hline
		\end{tabular}
	\end{table}

	\section{Conclusions and Outlook}\setcounter{equation}{0}
	
	The hybrid Fourier--finite difference scheme incorporates the magnetic Bloch boundary conditions through shifts of Fourier-mode indices and preserves the Hermitian structure of the discrete Hamiltonian. Under suitable regularity, spectral-gap, and discrete eigenspace dimension-preserving assumptions, the errors in the discrete eigenvalues and their associated eigenspaces are bounded by \(O(h_1^2+N_2^{-s})\). The numerical experiments show second-order convergence in the \(x\)-direction and rapid error reduction as the Fourier resolution increases. The tests include periodic approximations of multifrequency potentials and magnetic fields on enlarged magnetic unit cells, as well as computations of total densities for degenerate eigenspaces.
	
	Future work concerns the computation of Berry curvature and Chern numbers, together with associated error estimates, for isolated bands and multiple-eigenvalue eigenspaces separated from the rest of the spectrum by a gap. For magnetic fields or potentials containing incommensurate spatial frequencies, another direction is to develop structure-preserving discretizations and spectral approximations for the corresponding magnetic Schr\"odinger eigenvalue problems, drawing on higher-dimensional periodic representations and projection methods for quasiperiodic quantum systems \cite{Jiang-Li-Zhang2025,IWFPM2025}.


\begin{thebibliography}{99}
		
		\bibitem{AlougesBonnaillieNoel2006}
		F.~Alouges and V.~Bonnaillie-No\"el,
		Numerical computations of fundamental eigenstates for the Schr\"odinger operator under constant magnetic field,
		\emph{Numer. Methods Partial Differential Equations}, 22:1090--1105, 2006.
		doi: 10.1002/num.20137.
		
		\bibitem{AvronSimon1985}
		J.~E. Avron and B.~Simon,
		Stability of gaps for periodic potentials under variation of a magnetic field,
		\emph{J. Phys. A: Math. Gen.}, 18:2199--2205, 1985.
		doi: 10.1088/0305-4470/18/12/017.
		
		\bibitem{BabuskaOsborn1991}
		I.~Babu\v{s}ka and J.~Osborn,
		Eigenvalue problems,
		in P.~G. Ciarlet and J.~L. Lions (eds.), \emph{Handbook of Numerical Analysis, Vol.~II: Finite Element Methods (Part~1)},
		North-Holland, Amsterdam, 1991, pp.~641--787.
		doi: 10.1016/S1570-8659(05)80042-0.
		
		\bibitem{BonnaillieNoelDaugeMartinVial2007}
		V.~Bonnaillie-No\"el, M.~Dauge, D.~Martin and G.~Vial,
		Computations of the first eigenpairs for the Schr\"odinger operator with magnetic field,
		\emph{Comput. Methods Appl. Mech. Engrg.}, 196:3841--3858, 2007.
		doi: 10.1016/j.cma.2006.10.041.
		
		\bibitem{Brown1964}
		E.~Brown,
		Bloch electrons in a uniform magnetic field,
		\emph{Phys. Rev.}, 133:A1038--A1044, 1964.
		doi: 10.1103/PhysRev.133.A1038.
		
		\bibitem{CaiGalli2004}
		W.~Cai and G.~Galli,
		Ab initio calculations in a uniform magnetic field using periodic supercells,
		\emph{Phys. Rev. Lett.}, 92:186402, 2004.
		doi: 10.1103/PhysRevLett.92.186402.
		
		\bibitem{CancesEhrlacherMaday2014}
		E.~Canc\`es, V.~Ehrlacher and Y.~Maday,
		Non-consistent approximations of self-adjoint eigenproblems: Application to the supercell method,
		\emph{Numer. Math.}, 128:663--706, 2014.
		doi: 10.1007/s00211-014-0626-0.
		
		\bibitem{Christiansen}
		S.~H. Christiansen and T.~G. Halvorsen,
		A gauge invariant discretization on simplicial grids of the Schr\"odinger eigenvalue problem in an electromagnetic field,
		\emph{SIAM J. Numer. Anal.}, 49:331--345, 2011.
		doi: 10.1137/090757502.
		
		\bibitem{Cornean2010}
		H.~D. Cornean,
		On the Lipschitz continuity of spectral bands of Harper-like and magnetic Schr\"odinger operators,
		\emph{Ann. Henri Poincar\'e}, 11:973--990, 2010.
		doi: 10.1007/s00023-010-0048-1.
		
		\bibitem{DNR1978a}
		J.~Descloux, N.~Nassif and J.~Rappaz,
		On spectral approximation. Part~I: The problem of convergence,
		\emph{RAIRO Anal. Num\'er.}, 12:97--112, 1978.
		doi: 10.1051/m2an/1978120200971.
		
		\bibitem{FukuiHatsugaiSuzuki2005}
		T.~Fukui, Y.~Hatsugai and H.~Suzuki,
		Chern numbers in discretized Brillouin zone: Efficient method of computing (spin) Hall conductances,
		\emph{J. Phys. Soc. Jpn.}, 74:1674--1677, 2005.
		doi: 10.1143/JPSJ.74.1674.
		
		\bibitem{Gallistl2015}
		D.~Gallistl,
		An optimal adaptive FEM for eigenvalue clusters,
		\emph{Numer. Math.}, 130:467--496, 2015.
		doi: 10.1007/s00211-014-0671-8.
		
		\bibitem{GovernaleUngarelli1998}
		M.~Governale and C.~Ungarelli,
		Gauge-invariant grid discretization of the Schr\"odinger equation,
		\emph{Phys. Rev. B}, 58:7816--7821, 1998.
		doi: 10.1103/PhysRevB.58.7816.
		
		\bibitem{Gruber2000}
		M.~J. Gruber,
		Bloch theory and quantization of magnetic systems,
		\emph{J. Geom. Phys.}, 34:137--154, 2000.
		doi: 10.1016/S0393-0440(99)00059-5.
		
		\bibitem{HalvorsenKvaal2012}
		T.~G. Halvorsen and S.~Kvaal,
		Manifestly gauge invariant discretizations of the Schr\"odinger equation,
		\emph{Phys. Lett. A}, 376:1107--1114, 2012.
		doi: 10.1016/j.physleta.2012.02.028.
		
		\bibitem{HelfferSjostrand1989}
		B.~Helffer and J.~Sj\"ostrand,
		\'{E}quation de Schr\"odinger avec champ magn\'etique et \'equation de Harper,
		in H.~Holden and A.~Jensen (eds.), \emph{Schr\"odinger Operators}, Lecture Notes in Physics, Vol.~345,
		Springer, Berlin, 1989, pp.~118--197.
		doi: 10.1007/3-540-51783-9\_19.
		
		\bibitem{Hetmaniuk}
		U.~Hetmaniuk and R.~Lehoucq,
		Basis selection in LOBPCG,
		\emph{J. Comput. Phys.}, 218:324--332, 2006.
		doi: 10.1016/j.jcp.2006.02.007.
		
		\bibitem{Hofstadter1976}
		D.~R. Hofstadter,
		Energy levels and wave functions of Bloch electrons in rational and irrational magnetic fields,
		\emph{Phys. Rev. B}, 14:2239--2249, 1976.
		doi: 10.1103/PhysRevB.14.2239.
		
		\bibitem{Jiang-Li-Zhang2025}
		K. Jiang, S.-F. Li and P.-W. Zhang,
		On the approximation of quasiperiodic functions with Diophantine frequencies by periodic functions,
		\emph{SIAM J. Math. Anal.}, 57:951--978, 2025.
		doi: 10.1137/24M165925X.
		
		\bibitem{IWFPM2025}
		K. Jiang, X.-Y. Li, Y. Ma, J. Zhang, P.-W. Zhang and Q. Zhou,
		Irrational-Window-Filter Projection Method and Application to Quasiperiodic Schrödinger Eigenproblems.
		SIAM Journal on Numerical Analysis, 63(2), 564--587, 2025.
		DOI：10.1137/24M1666197
		
		\bibitem{Knyazev2001LOBPCG}
		A.~V. Knyazev,
		Toward the optimal preconditioned eigensolver: Locally optimal block preconditioned conjugate gradient method,
		\emph{SIAM J. Sci. Comput.}, 23:517--541, 2001.
		doi: 10.1137/S1064827500366124.
		
		\bibitem{LeeCaiGalli2007}
		E.~Lee, W.~Cai and G.~A. Galli,
		Electronic structure calculations in a uniform magnetic field using periodic supercells,
		\emph{J. Comput. Phys.}, 226:1310--1331, 2007.
		doi: 10.1016/j.jcp.2007.05.022.
		
		\bibitem{LiXiao2021}
		X.-Y.~Li and A.-G.~Xiao,
		Space-fractional diffusion equation with variable coefficients: Well-posedness and Fourier pseudospectral approximation,
		\emph{J. Sci. Comput.}, 87:28, 2021.
		doi: 10.1007/s10915-021-01435-w.
		
		\bibitem{Nenciu1991}
		G.~Nenciu,
		Dynamics of band electrons in electric and magnetic fields: Rigorous justification of the effective Hamiltonians,
		\emph{Rev. Mod. Phys.}, 63:91--127, 1991.
		doi: 10.1103/RevModPhys.63.91.
		
		\bibitem{OvallZhu2025}
		J.~S. Ovall and L.~Zhu,
		A canonical gauge for computing of eigenpairs of the magnetic Schr\"odinger operator,
		\emph{SIAM J. Sci. Comput.}, 47:B1026--B1044, 2025.
		doi: 10.1137/24M1692228.
		
		\bibitem{PanatiSpohnTeufel2003}
		G.~Panati, H.~Spohn, and S.~Teufel,
		Effective dynamics for Bloch electrons: Peierls substitution and beyond,
		\emph{Commun. Math. Phys.}, 242:547--578, 2003.
		doi: 10.1007/s00220-003-0950-1.
		
		\bibitem{Seye2025}
		A.~Seye and M.~Filoche,
		Localization structure of electronic states in the quantum Hall effect,
		\emph{Phys. Rev. B}, 112:144202, 2025.
		doi: 10.1103/266g-9zc3.
		
		\bibitem{ThoulessKohmotoNightingaleDenNijs1982}
		D.~J. Thouless, M.~Kohmoto, M.~P. Nightingale, and M.~den Nijs,
		Quantized Hall conductance in a two-dimensional periodic potential,
		\emph{Phys. Rev. Lett.}, 49:405--408, 1982.
		doi: 10.1103/PhysRevLett.49.405.
		
		\bibitem{Zak1964a}
		J.~Zak,
		Magnetic translation group,
		\emph{Phys. Rev.}, 134:A1602--A1606, 1964.
		doi: 10.1103/PhysRev.134.A1602.
		
	\end{thebibliography}
\end{document}